\documentclass{amsart}
\usepackage{amssymb}

\usepackage[v2]{xy}
\SilentMatrices

\usepackage{ifthen}
\usepackage{xspace}

\DeclareMathAlphabet{\scr}{U}{rsfs}{m}{n}

\newcommand{\op}{\mathrm{op}}

\newcommand{\Z}{\ensuremath{\mathbb{Z}}\xspace}

\newcommand{\ROG}{\ensuremath{RO(G)}\xspace}
\newcommand{\posROG}{\ensuremath{RO^{+}(G)}\xspace}
\newcommand{\tROG}{\widetilde{RO}(G)}
\newcommand{\free}[1]{\mathbb{#1}}
\newcommand{\cofree}[1]{\mathbb{#1}^{\sharp}}
\newcommand{\cell}[2]{R^{#1}[#2]}
\newcommand{\dcat}[1]{\scr{D}_{#1}}
\newcommand{\aC}{\scr{C}}

\newcommand{\HoGMod}{\mathrm{Ho}G\scr{M}}

\DeclareMathOperator{\Image}{Im}

\newcommand{\cofm}[1]{\bar\imath_{#1}}
\newcommand{\cofn}[1]{\bar\imath_{#1}}
\newcommand{\cofl}[1]{\bar\imath_{#1}}
\newcommand{\fibm}[1]{\bar\imath^{#1}}
\newcommand{\incm}[1]{f_{#1}}
\newcommand{\incn}[1]{f_{#1}}
\newcommand{\incl}[1]{f_{#1}}
\newcommand{\prom}[1]{f^{#1}}
\newcommand{\cofM}[1]{\bar\imath^{M}_{#1}}
\newcommand{\cofN}[1]{\bar\imath^{N}_{#1}}

\newcommand{\incM}[1]{f^{M}_{#1}}
\newcommand{\incN}[1]{f^{N}_{#1}}

\newcommand{\mac}[1]{\mathpalette\macunderline{#1}#1}
\newcommand{\macunderline}[2]{%
{\setbox0\hbox{$#1\mathrm{#2}$}\setbox2\hbox{$#1#2$}%
\ifdim\wd0<\wd2%
\hbox to 0pt{$\underline{\hbox to \wd0{\hss}}$\hss}%
\else%
\setbox2\hbox{$#1#2\!$}%
\hbox to 0pt{$\underline{\hbox to \wd2{\hss}}$\hss}%
\fi}}
\let\hmac\mac

\newcommand{\iso}{\cong}     
\newcommand{\sma}{\wedge}    

\newcommand{\from}{\leftarrow}
\newcommand{\overto}[1]{\xrightarrow{#1}}

\def\quickop#1{\expandafter\DeclareMathOperator\csname #1\endcsname{#1}}
\quickop{Hom} \quickop{Tel} \quickop{Mic}
\quickop{Ext} \quickop{Tor} \quickop{id} 
\quickop{Lim} \quickop{Colim} \quickop{Hocolim} \quickop{Holim}

\newcommand{\tindex}{\ensuremath{\tau}\xspace}

\newcommand{\aindex}{\ensuremath{\alpha}\xspace}
\newcommand{\bindex}{\ensuremath{\beta}\xspace}
\newcommand{\cindex}{\ensuremath{\gamma}\xspace}
\newcommand{\dindex}{\ensuremath{\delta}\xspace}

\newcommand{\shom}[3][{G}]{\ensuremath{{[ {#2},{#3} ]_{#1}}}\xspace}

\newcommand{\ePt}{\ensuremath{G/G}\xspace}

\newcommand{\mMBr}{\ensuremath{\mac{B}}\xspace}

\newcommand{\mML}{\ensuremath{\mac{L}}\xspace}
\newcommand{\mMM}{\ensuremath{\mac{M}}\xspace}
\newcommand{\mMN}{\ensuremath{\mac{N}}\xspace}
\newcommand{\mMP}{\ensuremath{\mac{P}}\xspace}

\newcommand{\mMR}{\ensuremath{\mac{R}}\xspace}

\newcommand{\mgMBr}[1][*]{\ensuremath{\mac{B}_{#1}}\xspace}

\newcommand{\mgMF}[1][*]{\ensuremath{\mac{F}_{#1}}\xspace}
\newcommand{\mgMI}[1][*]{\ensuremath{\mac{I}_{#1}}\xspace}
\newcommand{\mgMK}[1][*]{\ensuremath{\mac{K}_{#1}}\xspace}

\newcommand{\mgML}[1][*]{\ensuremath{\mac{L}_{#1}}\xspace}
\newcommand{\mgMM}[1][*]{\ensuremath{\mac{M}_{#1}}\xspace}
\newcommand{\mgMMp}[1][*]{\ensuremath{\mac{M}'_{#1}}\xspace}
\newcommand{\mgMMpp}[1][*]{\ensuremath{\mac{M}''_{#1}}\xspace}
\newcommand{\mgMN}[1][*]{\ensuremath{\mac{N}_{#1}}\xspace}
\newcommand{\mgMP}[1][*]{\ensuremath{\mac{P}_{#1}}\xspace}
\newcommand{\mgMQ}[1][*]{\ensuremath{\mac{Q}_{#1}}\xspace}
\newcommand{\mgMR}[1][*]{\ensuremath{\mac{R}_{#1}}\xspace}

\let\pF\mgMF
\let\pI\mgMI
\let\pK\mgMK
\let\pL\mgML
\let\pM\mgMM
\let\pN\mgMN

\let\pP\mgMP

\let\pR\mgMR

\newcommand{\mzgMM}[1][\#]{\mgMM[#1]}
\newcommand{\mzgMN}[1][\#]{\mgMN[#1]}

\newcommand{\mbMC}[2]{\ensuremath{\mac{C}_{{#1},{#2}}}\xspace}
\newcommand{\mbMD}[2]{\ensuremath{\mac{D}_{{#1},{#2}}}\xspace}
\newcommand{\mbMO}[2]{\ensuremath{\mac{O}_{{#1},{#2}}}\xspace}
\newcommand{\mbMP}[2]{\ensuremath{\mac{P}_{{#1},{#2}}}\xspace}
\newcommand{\mbMQ}[2]{\ensuremath{\mac{Q}_{{#1},{#2}}}\xspace}

\newcommand{\mcMC}[2]{\ensuremath{\mac{C}^{#1}_{#2}}\xspace}
\newcommand{\mcMD}[2]{\ensuremath{\mac{D}^{#1}_{#2}}\xspace}
\newcommand{\mcMI}[2]{\ensuremath{\mac{I}^{#1}_{#2}}\xspace}

\newcommand{\pcR}[1][*]{\hmac{R}^{#1}}
\newcommand{\pcM}[1][*]{\hmac{M}^{#1}}

\newcommand{\gsL}{\ensuremath{L}\xspace}
\newcommand{\gsM}{\ensuremath{M}\xspace}
\newcommand{\gsN}{\ensuremath{N}\xspace}
\newcommand{\gsS}{\ensuremath{S}\xspace}

\newcommand{\Burn}{\ensuremath{\mathfrak{B}_{G}}\xspace}
\newcommand{\OBurn}{\ensuremath{\mathfrak{O}_{G}}\xspace}
\newcommand{\AB}{\ensuremath{\mathfrak{Ab}}\xspace}
\newcommand{\Mack}{\ensuremath{\mathfrak{M}}\xspace}
\newcommand{\RMack}{\ensuremath{\mathfrak{M}_*}\xspace}
\newcommand{\ZMack}{\ensuremath{\mathfrak{M}_{\#}}\xspace}
\newcommand{\rMod}[1][{\mgMR}]{\ensuremath{{{#1}\textnormal{-Mod}}}\xspace}
\newcommand{\Modr}[1][{\mgMR}]{\ensuremath{{\textnormal{Mod-}{#1}}}\xspace}

\newcommand{\inj}[2]{\ensuremath{\mac{\mathcal{I}}({#1},{#2})}\xspace}
\newcommand{\ginj}[3][*]{\ensuremath{\mac{\mathcal{I}}({#2},{#3})_{#1}}\xspace}
\newcommand{\func}[2]{\ensuremath{\langle {#1},{#2}\rangle}\xspace}

\newcommand{\pimac}[1][*]{\ensuremath{\mac{\pi}_{#1}}\xspace}
\newcommand{\pimacS}[2][*]{\ensuremath{\pimac[#1]({#2})}\xspace}

\newcommand{\boxprod}{\mathbin{\mathpalette{}\square}}
\newcommand{\gbxprd}{\ensuremath{\boxprod_{*}}}
\newcommand{\zgbxprd}{\ensuremath{\boxprod_{\#}}}
\newcommand{\rbox}[1][{\mgMR}]{\ensuremath{\boxprod_{#1}}}

\newcommand{\ghom}[3][*]{\ensuremath{{\langle {#2},{#3} \rangle_{#1}}}\xspace}
\newcommand{\rfunc}[4][{}]{\ensuremath{{\langle {#3},{#4} \rangle^{{#2}}_{#1}}}\xspace}

\def\dersma{\Tor}
\def\derF{\Ext}

\DeclareMathOperator{\TTor}{\mskip.5\thinmuskip\hmac{\mskip-.5\thinmuskip Tor\mskip-.5\thinmuskip}\mskip.5\thinmuskip}
\DeclareMathOperator{\TExt}{\mskip.25\thinmuskip\hmac{\mskip-.25\thinmuskip Ext\mskip-.125\thinmuskip}\mskip.125\thinmuskip}

\let\MTor\TTor
\let\MExt\TExt
\newcommand{\mExt}[4][*]{\ensuremath{\MExt^{#1}_{#2}({#3},{#4})}\xspace}
\newcommand{\mTor}[4][*]{\ensuremath{\MTor^{#2}_{#1}({#3},{#4})}\xspace}

\newcommand{\cocycle}[1]{\ensuremath{\mathfrak{a}_{#1}}\xspace}
\newcommand{\dclass}[1]{\ensuremath{\mathfrak{b}_{#1}}\xspace}
\newcommand{\INV}{\mathfrak{S}}
\newcommand{\SPH}{\mathfrak{C}}
\newcommand{\Bunits}{A^{*}}
\newcommand{\coclass}{\bar{\mathfrak{a}}}
\newcommand{\poscoclass}{\bar{\mathfrak{a}}^{+}}

\newlength{\arsize}

\newcommand{\lrarrow}[1][{\:}]{\settowidth{\arsize}{\ensuremath{\;\;\;\;{{}^{#1}}}}%
	\ensuremath{\:{\xy \ar @{->}^-{#1}<\arsize,0em> \endxy}\;\,}}

\newcommand{\llarrow}[1][{\;}]{\settowidth{\arsize}{\ensuremath{\;\;\;\;{{}^{#1}}}}%
	\ensuremath{\:{\xy \ar @{<-}^<>(.65){#1}<\arsize,0em> \endxy}\;\;}}
\newcommand{\Bmap}[4][{\lrarrow}]{\ensuremath{{#2} : {#3}{#1}{#4}}}
\newcommand{\map}[3]{\Bmap{#1}{#2}{#3}}

\makeatletter
\newcounter{thmpart}

\newenvironment{thmlist}[1][0]%
{\begin{list}%
{\ifthenelse{\boolean{@newlist} \and \(\not {\equal{#1}{0}}\)}%
	{\hspace*{-2em} \normalfont \the\thm@headfont \upshape (\alph{thmpart})}%
	{\normalfont \the\thm@headfont \upshape (\alph{thmpart})}%
}%
{\usecounter{thmpart}%
\setlength{\labelsep}{0.5em}\setlength{\leftmargin}{0pt}%
\setlength{\labelwidth}{0pt}\setlength{\itemindent}{1.8em}%
\setlength{\listparindent}{1em}}}%
{\end{list}}

\newcounter{hyppart}

{\begin{list}%
{\ifthenelse{\boolean{@newlist} \and \(\not {\equal{#1}{0}}\)}%
	{\hspace*{-2em} \normalfont \the\thm@headfont \upshape (\roman{hyppart})}%
	{\normalfont \the\thm@headfont \upshape (\roman{hyppart})}%
}%
{\usecounter{hyppart}%
\setlength{\labelsep}{0.5em}\setlength{\leftmargin}{0pt}%
\setlength{\labelwidth}{0pt}\setlength{\itemindent}{1.8em}%
\setlength{\listparindent}{1em}}}%
{\end{list}}

\def\llabel#1{\@bsphack
  \protected@write\@auxout{}%
         {\string\newlabel{#1}{{\alph{\@listctr}}{\thepage}}}%
  \@esphack}

\def\hlabel#1{\@bsphack
  \protected@write\@auxout{}%
         {\string\newlabel{#1}{{\roman{\@listctr}}{\thepage}}}%
  \@esphack}
\makeatother

\newtheorem*{theorem}{Theorem}

\newtheorem{thm}{Theorem}[section]
\newtheorem*{thm*}{Theorem}
\newtheorem{cor}[thm]{Corollary}
\newtheorem{lem}[thm]{Lemma}
\newtheorem{prop}[thm]{Proposition}

\theoremstyle{definition}
\newtheorem{defn}[thm]{Definition}

\theoremstyle{remark}
\newtheorem*{remark}{Remark}
\newtheorem{rem}[thm]{Remark}

\numberwithin{equation}{section}

\makeatletter
\let\c@equation\c@thm
\makeatother

\begin{document}
\title[Equivariant Spectral Sequences]
{Equivariant Universal Coefficient and K\"unneth Spectral Sequences}

\author{L. Gaunce Lewis, Jr.}
\address{Department of Mathematics, Syracuse University, 
Syracuse, NY \ 13244-1150}
\email{lglewis@syr.edu}

\author{Michael A. Mandell}
\address{DPMMS, University of Cambridge, Wilberforce Road, Cambridge CB3 0WB, UK}
\email{M.A.Mandell@dpmms.cam.ac.uk}
\thanks{The second author was supported in part by NSF grant DMS-0203980}

\date{October 4, 2004}

\begin{abstract}
We construct hyper-homology spectral sequences of \Z-graded and
\ROG-graded Mackey functors for $\Ext$ and $\Tor$ over $G$-equivariant
$S$-algebras ($A_{\infty}$ ring spectra) for finite groups $G$.  These
specialize to universal coefficient and K\"unneth spectral sequences.
\end{abstract}

\maketitle


\section*{Introduction}

In the non-equivariant context, universal coefficient and K\"unneth
spectral sequences provide important tools for computing generalized
homology and cohomology.  EKMM~\cite{EKMM} constructs examples of
these types of spectral sequences for theories that come from
``$S$-algebras'' (or, equivalently, $A_{\infty}$ ring spectra).  These
spectral sequences are special cases of the general ``hyper-homology''
spectral sequences for computing ``$\Tor$'' (the homotopy groups of the
derived smash product) and ``$\Ext$'' (the homotopy groups of the
derived function module) for the category of modules over an 
$S$-algebra.  The purpose of this
paper is to construct versions of these hyper-homology spectral 
sequences for the category of modules over an 
$G$-equivariant $S$-algebra
indexed on a complete universe, where $G$ is a finite group.  In
Section~\ref{secuct}, we derive a
number of equivariant universal coefficient and K\"unneth spectral sequences from
these equivariant hyper-homology spectral sequences.

In this context, the homotopy groups of a $G$-spectrum $X$ form a
graded Mackey functor.  There is an abelian group 
$(\pimac[q]X)(G/H)=\pi_{q}^{H}X=\pi_{q}(X^{H})$ for each subgroup $H$
of $G$, and for each subgroup $K$ of $H$ there are maps $\pi_{q}^{H}X\to
\pi_{q}^{K}X$ (induced by the inclusion of fixed points) and
$\pi_{q}^{K}X\to \pi_{q}^{H}X$ (the transfer), satisfying various
natural relations.  This ``homotopy
Mackey functor'' can be regarded as graded over either the integers or the
real representation ring of $G$.  Both choices lead to equivariant
generalizations of the spectral sequences of EKMM~\cite{EKMM}.  
Although these generalizations appear formally similar and are developed 
here in parallel, the spectral sequences have a very different
calculational feel.

As reviewed in Section~\ref{secgrdmackey}, the categories of 
\Z-graded and of \ROG-graded Mackey functors have
symmetric monoidal products denoted ``$\zgbxprd$'' and ``$\gbxprd$''
respectively.  The products have adjoint function object functors
``$\ghom[\#]{-}{-}$'' and ``$\ghom{-}{-}$''.  There is an obvious 
notion of a graded Mackey functor ring (in either context), which
consists of a graded Mackey functor $\mgMR$ with an associative and
unital multiplication $\mgMR\boxprod \mgMR\to \mgMR$.  We have the
evident notion of a left (resp. right) $\mgMR$-module as a graded
Mackey functor $\mgMM$ with an associative and unital map
$\mgMR\boxprod \mgMM\to \mgMM$ (resp. $\mgMN\boxprod \mgMR\to \mgMN$).
The usual coequalizer defines a functor taking a 
right $\mgMR$-module $\mgMN$ and a left $\mgMR$-module $\mgMM$ to 
a graded Mackey functor $\mgMN\rbox[\mgMR] \mgMM$.  Similarly, the usual
equalizer defines a functor taking a pair of left
$\mgMR$-modules \mgML and \mgMM to the graded Mackey functor $\rfunc{\mgMR}{\mgML}{\mgMM}.$  In Section~\ref{sechom},
standard homological algebra is used to construct the associated 
derived functors.  We 
define $\MTor^{\mgMR}_{s}(\mgMN,\mgMM)$ as the $s$-th left derived
functor of $\mgMN\rbox[\mgMR]\mgMM$ and
$\MExt_{\mgMR}^{s}(\mgML,\mgMM)$ as the $s$-th right derived functor
of $\rfunc{\mgMR}{\mgML}{\mgMM}$.  These are regarded as bigraded
Mackey functors with the usual conventions:
$\MTor^{\mgMR}_{s,\tindex}(\mgMN,\mgMM) =(\MTor^{\mgMR}_{s}(\mgMN,\mgMM))_{\tindex}$ and
$\MExt_{\mgMR}^{s,\tindex}(\mgML,\mgMM) =
(\MExt_{\mgMR}^{s}(\mgML,\mgMM))_{-\tindex}$. The homological grading
in $s$ is always over the non-negative integers whereas the internal
grading in $\tindex$ is over \Z\ or \ROG.  When $\mgMR$ is commutative, 
$\rbox[\mgMR]$ becomes a closed symmetric monoidal product on the
category of left $\mgMR$-modules, and $\MTor^{\mgMR}_{s}$ and
$\MExt_{\mgMR}^{s}$ become graded $\mgMR$-modules.  More generally,
when $M$ is an $\mgMR$-bimodule (has commuting left and right module
structures), $\MTor^{\mgMR}_{s}(\mgMN,\mgMM)$ and
$\MExt_{\mgMR}^{s}(\mgML,\mgMM)$ become graded right $\mgMR$-modules
via the ``unused'' $\mgMR$-module structure on $\mgMM$.  Similarly,
when $\mgMN$ and $\mgML$ are $\mgMR$-bimodules,
$\MTor^{\mgMR}_{s}(\mgMN,\mgMM)$ and $\MExt_{\mgMR}^{s}(\mgML,\mgMM)$
become graded left $\mgMR$-modules.

The equivariant stable category has a symmetric monoidal
product ``$\sma$'', and the homotopy Mackey functor $\pimac $ is a
lax symmetric monoidal functor into either \Z-graded or \ROG-graded 
Mackey functors.  In other words,
there is a suitably associative, commutative, and unital natural
transformation 
\[
\pimac X \boxprod \pimac Y \lrarrow \pimac (X\sma Y).
\]
For formal reasons, we obtain an adjoint natural transformation
\[
\pimac F(X,Y)\lrarrow \func{\pimac X}{\pimac Y},
\]
where $F$ denotes the function spectrum.  
Further, when $R$ is a homotopical ring spectrum (i.e., a ring
spectrum in the equivariant stable category),
$\pimac R$ is a Mackey functor ring, commutative when $R$ is.
Likewise, when $M$ is a homotopical $R$-module, $\pimac M$ becomes
a $\pimac R$-module.  When $R$ is an equivariant $S$-algebra, we have a
smash product over $R$ and an $R$-function spectrum functor.  In this
context, it is convenient to extend the conventions of EKMM
\cite{EKMM} by using $\TTor^{R}_{*}(N,M)$ and $\TExt_{R}^{*}(L,M)$ for the
homotopy Mackey functors of the derived
functors.  The natural transformations above descend to natural
transformations 
\[
\pN \rbox[\pR] \pM \to \TTor^{R}_{*}(N,M)
\quad\text{and}\quad
\TExt_{R}^{*}(L,M)\to \rfunc{\pR}{\pL}{\pM},
\]
where $\pR=\pimac R$, $\pL=\pimac L$, $\pM=\pimac M$, and
$\pN=\pimac N$.  When $M$ is a weak bimodule, i.e., has a homotopical
right $R$-module structure in the category of left $R$-modules,
$\TTor^{R}_{*}(N,M)$ and $\TExt_{R}^{-*}(L,M)$ are naturally right
$\pR$-modules and the maps above are maps of right $R$-modules.  
Similar statements hold for $L$ and $N$.
Our main results are the following two theorems.

\begin{theorem}(Hyper-Tor Spectral Sequences)
Let $G$ be a finite group, $R$ be an equivariant $S$-algebra indexed on a
complete universe, and $M$ and $N$ be a left and a right 
$R$-module.  There are strongly convergent spectral sequences  
\[
E^{2}_{s,\tindex}=\MTor^{\pR}_{s,\tindex}(\pN,\pM) 
\quad \Longrightarrow \quad 
\TTor^{R}_{s+\tindex}(N,M),
\]
of \Z-graded Mackey functors and of \ROG-graded Mackey functors which are 
natural in $M$ and $N$.  The edge homomorphisms are the usual natural
transformations $\pN \rbox[\pR] \pM\to \TTor^{R}_{*}(N,M)$.
When $R$ is commutative or $M$ is a weak bimodule, these
are spectral sequences of right $\pR$-modules. When $N$ is a weak
bimodule, they are spectral sequences of left $\pR$-modules.
\end{theorem}

\begin{theorem}(Hyper-Ext Spectral Sequences) Let $G$ be a finite
group, $R$ be an equivariant $S$-algebra indexed on a complete
universe, and $L$ and $M$ be left $R$-modules.  There are 
conditionally convergent spectral sequences 
\[
E_{2}^{s,\tindex}=\MExt_{\pR}^{s,\tindex}(\pL,\pM) 
\quad \Longrightarrow \quad 
\TExt_{R}^{s+\tindex}(L,M),
\]
of \Z-graded Mackey functors and of \ROG-graded Mackey functors which are 
natural in $L$ and $M$.  The edge homomorphisms are the usual natural
transformations $\TExt_{R}^{*}(L,M)\to \rfunc{\pR}\pL\pM$.
When $R$ is commutative or $M$ is a weak bimodule,
these are spectral sequences of right $\pcR$-modules.  When $L$ is a weak
bimodule, these are spectral sequences of left $\pcR$-modules.
\end{theorem}


To be precise about the differentials, the Hyper-Tor spectral sequence
has $r$-th differential 
\[
d^{r}_{s,t}\colon E^{r}_{s,t}\lrarrow E^{r}_{s-r,t+r-1}.
\]
Thus, in the \Z-graded context, it is a homologically graded right
half-plane spectral sequence of Mackey functors.  In the \ROG-graded
context, it essentially consists of one homologically graded right half-plane
spectral sequence of Mackey functors for each element of the subgroup 
$\tROG$ of \ROG consisting of  virtual representations of dimension zero.  The
Hyper-Ext spectral sequence has $r$-th differential 
\[
d_{r}^{s,t}\colon E_{r}^{s,t}\lrarrow E_{r}^{s+r,t-r+1}.
\]
In the \Z-graded context, it is a cohomologically graded right 
half-plane spectral sequence of Mackey functors.  In the \ROG-graded
context, it consists of one cohomologically graded right half-plane
spectral sequence of Mackey functors for each element of $\tROG$.

In the \Z-graded context, ``strong'' convergence means that for 
every subgroup $H$ of $G$, the spectral
sequence of abelian groups $E^{r}_{s,t}(G/H)$ 
converges strongly.
In the \ROG-graded context, it means that for every subgroup $H$ of
$G$ and every $\tilde\tindex$ in $\tROG$, the spectral sequence consisting of the abelian groups 
$E^{r}_{s,\tilde\tindex+t}(G/H)$ 
converges strongly.  Conditional convergence is defined analogously in terms of the abelian groups $E_{r}^{s,t}(G/H)$ and $E_{r}^{s,\tilde\tindex+t}(G/H)$.  

We emphasize again that, although we construct them in parallel, the
\Z-graded and \ROG-graded versions of these spectral sequences are 
typically very different.  For example, if $L=R\sma S^{V}$ for a
finite dimensional $G$-representation $V$, then $\pL$ is projective as
an \ROG-graded Mackey functor $\pR$-module, but is usually not
projective as a \Z-graded Mackey functor 
$\pR$-module.  The \ROG-graded spectral sequence 
collapses at $E_{2}$ and is concentrated in homological degree zero;
this typically does not happen in the \Z-graded spectral sequence.

The peculiarities in the homological algebra of Mackey
functors that occur when $G$ is not finite or when the universe is not
complete (cf. Lewis \cite{PrjFlt}) make it worthwhile to observe that 
the homological behavior of the graded Mackey functors in the present
context is perfectly standard.  

\begin{remark}
Let $\mgMR$ be a \Z-graded or \ROG-graded Mackey functor ring.  
\begin{thmlist}
\item If $\mgMM$ is a projective left $\mgMR$-module or $\mgMN$ is a
projective right $\mgMR$-module, then $\MTor^{\mgMR}_{s}(\mgMN,\mgMM)=0$ in
homological degree $s>0$. In this case, the Hyper-Tor spectral
sequence above collapses at $E^{2}$ and the edge homomorphism is an
isomorphism. 

\item If $\mgML$ is a projective left $\mgMR$-module or $\mgMM$
is an injective left $\mgMR$-module, then $\MExt^{s}_{\mgMR}(\mgML,\mgMM)=0$
in cohomological degree $s>0$.  In this case, the Hyper-Ext spectral
sequence above collapses at $E_{2}$ and the edge homomorphism is an
isomorphism. 

\item For any left $\mgMR$-modules $\mgML$, $\mgMM$, the usual (Yoneda)
$\Ext$ groups can be identified as
$\Ext^{s}_{\mgMR}(\mgML,\mgMM)=(\MExt^{s,0}_{\mgMR}(\mgML,\mgMM))(G/G)$.
Evaluation of a Mackey functor at $G/G$ is exact.  Thus, our Hyper-Ext
spectral sequence specializes to a conditionally convergent spectral
sequence
\[
E_{2}^{s,t}=\Ext^{s,t}_{\pR}(\pL,\pM)
\quad \Longrightarrow \quad 
\HoGMod^{s+t}_{R}(L,M).
\]
Here, $\HoGMod^{s+t}_{R}(L,M)$ is the abelian group of 
$G$-maps from $L$ to
$\Sigma^{s+t}M$ in the homotopy category of $R$-modules. 
\end{thmlist}
\end{remark}

The final topic considered in this paper is that of composition 
pairings.  The internal function
object $F_{R}(-,-)$ has a composition pairing
\[
F_{R}(M,N) \sma F_{R}(L,M)\lrarrow F_{R}(L,N).
\]
that induces a pairing
\[
\TExt^{s}_{R}(M,N) \boxprod \TExt^{t}_{R}(L,M)\lrarrow \TExt^{s+t}_{R}(L,N)
\]
If $R$ is commutative, then this composition pairing descends to
$\sma_{R}$, becoming a map of $R$-modules, and the pairing on $\TExt$
descends to $\rbox[\pcR]$, becoming a map of $\pcR$-modules.  Likewise,
for graded Mackey functor modules, the internal function object
$\rfunc{\mgMR}{-}{-}$ has a composition pairing that induces a
``Yoneda pairing'' on $\MExt$-objects
\[
\MExt^{s}_{\mgMR}(\mgMM,\mgMN)\boxprod \MExt^{t}_{\mgMR}(\mgML,\mgMM)
\lrarrow
\MExt^{s+t}_{\mgMR}(\mgML,\mgMN).
\]
We prove the following theorem regarding these pairings. 

\begin{theorem}
The Yoneda pairing on $\MExt^{*}_{\pR}$ induces a pairing of Hyper-Ext
spectral sequences that 
converges (conditionally) to the composition pairing on
$\TExt^{*}_{R}$.  When $R$ is commutative, this is a pairing of
spectral sequences of $\pcR$-modules.
\end{theorem}

In the case when $G$ is the trivial group, our argument corrects an error
in EKMM \cite[IV\S5]{EKMM}.

\subsection*{Organization}
In section~\ref{secuct}, we construct various universal coefficient and K\"unneth
spectral sequences from the spectral sequences described above.
Section~\ref{secgrdmackey} contains a brief review of the categories of
graded Mackey functors, graded box products, and graded function
Mackey functors.  The definitions of graded Mackey functor rings and modules and 
the definitions and basic properties of the graded box product and
function object over a graded Mackey functor ring are reviewed in Section~\ref{secgrdmods}.  
In Section~\ref{sechom}, we discuss projective and injective graded 
Mackey functor modules and use them to define $\MExt$ and $\MTor$.  Our work in
equivariant stable homotopy theory begins in Section~\ref{sectop} with a 
discussion of $R$-modules whose
homotopy Mackey functors are projective or injective $\pR$-modules.  
The spectral sequences introduced above are constructed in Section~\ref{secss}.  
Our results on the naturality of these spectral sequences and on the Yoneda 
pairing are proven in Section~\ref{secyon}.  Finally, in the appendix, we 
prove the folk theorem that the \ROG-graded
homotopy Mackey functor is a lax symmetric monoidal functor, i.e., 
that the smash product of equivariant spectra is compatible with the
graded box product.

\subsection*{Notation and conventions} Throughout this paper, $G$ is a
fixed finite group, and the ``equivariant stable category'' means the
derived category (obtained by formally inverting the weak
equivalences) of one of the modern categories of equivariant spectra,
i.e., the category of equivariant $S$-modules, the category of
equivariant orthogonal spectra, or the category of equivariant
symmetric spectra.  Except in Section~\ref{secuct}, which is written
from the point of view of homology and cohomology theory, the indexing
universe for spectra is always understood to be complete.  The term
``$G$-spectrum'' is used as an abbreviation for ``object in the
equivariant stable category''.

Except in Section~\ref{secyon}, we phrase all statements and arguments in
the equivariant stable category and in derived categories of modules.
The undecorated symbol ``$\sma$'' means the smash
product of spectra in the equivariant stable category or the smash
product with a spectrum in the derived category of $R$-modules.  We
write ``$\sma_{S}$'' when we mean the point-set smash product.
Likewise $F(-,-)$ denotes the derived function spectrum.  The point-set 
level functor adjoint to $\sma_{S}$ is denoted $F_{S}(-,-)$.  If $R$ is 
an equivariant $S$-algebra, then $\dersma^{R}(N,M)$ is the
derived balanced smash product of a left $R$-module $M$ and a right
$R$-module $N$.  Also, $\derF_{R}(L,M)$ is the derived function
spectrum (or $R$-module) of left $R$-modules $L$ and $M$.  Further, 
$\TTor^{R}_{*}(M,N) = \pimac (\dersma^{R}(N,M))$ and $\TExt_{R}^{*} =
\pimac[-*](\derF_{R}(L,M))$.  We reserve $(-)\sma_{R}(-)$ and
$F_{R}(-,-)$ for the point-set level functors.

\subsection*{Grading conventions}
The usual conventions for homology and cohomology theories and
for derived functors require the introduction of ``homological'' and
``cohomological'' grading.  For a (homologically) graded Mackey
functor $\mgMM$, the
corresponding cohomologically graded Mackey functor $\mMM^{*}$ is 
defined by $\mMM^{\aindex}=\mgMM[-\aindex]$ for all $\alpha$.  When $\mgMM$ 
is a left $\mgMR$-module, $\mMM^{*}$ is a left $\mMR^{*}$-module.
The choice of whether to regard the category of left $\mgMR$-modules and the
category of left $\mMR^{*}$-modules as different categories or the
same category with two different grading conventions is a
philosophical one: 
The same formulas hold provided we are consistent about
which grading we use in the latter case.  For example, there are
canonical natural isomorphisms 
\begin{gather*}
\MExt^{s,\tindex}_{\mMR^{*}}(\mML^{*},\mMM^{*})
\iso \MExt^{s,\tindex}_{\mgMR}(\mgML,\mgMM)
\\
\MTor_{s,\tindex}^{\mMR^{*}}(\mMN^{*},\mMM^{*})
\iso \MTor_{s,\tindex }^{\mgMR}(\mgMN,\mgMM)
\end{gather*}
(where $\MTor_{s,\tindex}^{\mMR^{*}}(\mMN^{*},\mMM^{*})=
(\MTor_{s}^{\mMR^{*}}(\mMN^{*},\mMM^{*}))^{-\tindex}$).
For this reason, except in the statements of the spectral sequences in
Section~\ref{secuct} and in the final construction of them in
Section~\ref{secss}, the main body of the paper treats graded
Mackey functors exclusively in homologically graded terms.


\section{Universal Coefficient and K\"unneth Spectral Sequences}
\label{secuct}

The hyper-homology spectral sequences of the introduction lead formally
to spectral sequences for computing equivariant homology and
cohomology.  In
particular, we obtain universal coefficient spectral sequences for 
computing $\hmac{E}_{*}(X)$ and $\hmac{E}^{*}(X)$ when either the theory
$E$ is represented by an $R$-module or the spectrum $X$ is represented
by an $R$-module.  Likewise, there are K\"unneth spectral sequences for computing 
$\hmac{R}_{*}(X\sma Y)$ and $\hmac{R}^{*}(X\sma Y)$.  Since the
equivariant homology and cohomology theories involved are \ROG-graded,
the arguments apply to the construction of both 
\ROG-graded and \Z-graded spectral sequences.  The resulting 
\ROG-graded spectral sequences are new even in the case when the
theories involved are ordinary (Bredon) homology and cohomology.  In
fact, the motivation for the research leading to this paper was the
need for such spectral sequences in \cite{ZpLinFree}.  The K\"unneth
spectral sequences are new for equivariant $K$-theory.

In the statements below, $R$ is a fixed
equivariant $S$-algebra indexed on a complete universe.  Given a
$G$-spectrum $E$ indexed on a complete universe and a $G$-spectrum
$X$ indexed on any universe $U'$, the $E$-homology of $X$ is the
\ROG-graded homotopy Mackey functor of $(i_{*}X)\sma E$.  Here, $i_{*}$ is
the left adjoint of the forgetful functor from the equivariant stable
category indexed on the complete universe to the equivariant stable
category indexed on $U'$.  The spectrum $(i_{*}X)\sma E$ is sometimes denoted  
$X\sma E$ to compactify notation.  The $E$-cohomology of $X$ is the
\ROG-graded Mackey functor
\[
\hmac{E}^{\alpha}(X)(G/H) 
= \shom{X \sma G/H_{+}}{\Sigma^{\alpha}E}
\iso \shom[H]{X}{\Sigma^{\alpha}E}
\]
or more accurately $\shom{i_{*}X \sma G/H_{+}}{\Sigma^{\alpha}E}$.  Of
course, all functors used here should be understood as the derived
functors. 

The forgetful functor from the derived category of left $R$-modules to
the equivariant stable category on a complete universe has a left
adjoint functor that we denote as $\free{R}$.  We write
$\free{R}^{\op}$ for the analogous functor into the derived category of
right $R$-modules.  These functors have the
usual properties expected of derived free $R$-module functors
(q.v. \cite{LMmmmc}).  For example, since $R$-modules are 
necessarily homotopical $R$-modules, we have natural maps 
\[
R\sma X \to \free{R}X
\qquad \text{and}\qquad
X\sma R \to \free{R}^{\op}X
\]
in the equivariant stable category (on a complete
universe).  These maps are always isomorphisms in the equivariant 
stable category.  We can compose the unit
map $X\to \free{R}X$ with the canonical comparison maps for the
derived smash product and derived function spectra to get natural maps 
\[
X \sma M \to \dersma^{R}({\free{R}^{\op}X},M)
\quad \text{and}\quad 
\derF_{R}({\free{R}X},M)\to F(X,M)
\]
in the equivariant stable category for any left $R$-module $M$.  These
maps are also always isomorphisms.  Just as $R\sma X$ and $X\sma R$
are naturally homotopical $R$-bimodules ($R$-bimodules in the
equivariant stable category), $\free{R}X$ and $\free{R}^{\op}X$ are
naturally weak bimodules in their respective categories: $\free{R}X$
is a homotopical right $R$-module in the derived category of left
$R$-modules and $\free{R}^{\op}X$ is a homotopical left $R$-module in
the derived category of right $R$-modules.  Moreover, the four maps in the
equivariant stable category displayed above are maps of homotopical
$R$-modules.  

As a consequence, by taking $N = \free{R}^{\op}X$ in the
Hyper-Tor spectral sequence and $L = \free{R}X$ in the Hyper-Ext 
spectral sequence, we obtain the following
spectral sequences of $\pR$-modules.

\begin{thm}(Universal Coefficient Spectral Sequence)
\label{uct}
Let $X$ be a $G$-spectrum and $M$ be a left $R$-module.
There is a natural strongly convergent homology spectral sequence of 
$\pR$-modules 
\[
E^{2}_{s,\tindex}=\MTor^{\pR}_{s,\tindex}(\hmac{R}_{*} X,\pM)
\quad \Longrightarrow \quad 
\hmac{M}_{s+\tindex}X
\]
and a natural conditionally convergent cohomology
spectral sequence of $\pR$-modules
\[
E_{2}^{s,\tindex}=\MExt_{\pR}^{s,\tindex}(\hmac{R}_{*}X,\pM) 
\quad \Longrightarrow \quad 
\hmac{M}^{s+\tindex}X.
\]
\end{thm}

The forgetful functor from the derived category of left $R$-modules to
the equivariant stable category also has a right adjoint functor. This 
adjoint is denoted $\cofree{R}$ and has the usual
properties expected from a derived cofree $R$-module functor
\cite{LMmmmc}.  In particular,  
if $E$ is a $G$-spectrum indexed on a complete universe, then the canonical map
$\cofree{R}E\to F(R,E)$ and the canonical comparison map of derived
function spectra 
\[
\derF_{R}(M,{\cofree{R}E}) \lrarrow F(M,E)
\]
are isomorphisms in the equivariant stable category.  Once again,
$\cofree{R}E$ is naturally a weak bimodule and the above maps in the
equivariant stable category are maps
of homotopical $R$-modules.
Plugging
$\cofree{R}E$ into the Hyper-Ext spectral sequence then gives the
cohomological spectral sequence in the theorem below.  Under the natural 
isomorphism 
\[
\hmac{R}_{*}E = \pimac  (E\sma R) 
\iso \pimac  (R\sma E) = \hmac{E}_{*}R,
\]
the homological spectral sequence in this theorem coincides
with the homological spectral sequence of the previous theorem.

\begin{thm}\label{auct}
Let $E$ be a $G$-spectrum
indexed on a complete universe and 
$M$ be a left $R$-module.
There is a natural strongly convergent homology spectral sequence of 
$\pR$-modules
\[
E^{2}_{s,\tindex}=\MTor^{\pR}_{s,\tindex}(\hmac{E}_{*} R,\pM) 
\quad \Longrightarrow \quad 
\hmac{E}_{s+\tindex}M
\]
and a natural conditionally convergent cohomology spectral sequence of $\pR$-modules
\[
E_{2}^{s,\tindex}=\MExt_{\pcR}^{s,\tindex}(\pcM, \hmac{E}^{*}R) 
\quad \Longrightarrow \quad 
\hmac{E}^{s+\tindex}M.
\]
\end{thm}

K\"unneth spectral sequences for computing the homology and
cohomology of $X\sma Y$ arise as special cases of our universal coefficient
spectral sequences.  Let $X$ and $Y$ be $G$-spectra indexed on the
same universe $U'$.  The change of universe
functor $i_{*}$ is symmetric monoidal; that is, there is an
isomorphism $i_{*}X\sma i_{*}Y\iso
i_{*}(X\sma Y)$ in the equivariant stable category.  The 
observations preceding Theorem~\ref{uct} yield a canonical isomorphism 
\[
i_{*}(X\sma Y) \sma R \iso (i_{*}X)\sma (i_{*}Y)\sma R
\iso (i_{*}X)\sma R\sma (i_{*}Y)
\iso \dersma^{R}(\free{R}^{\op}(i^{*}X),\free{R}(i^{*}Y))
\]
in the equivariant stable category.
Taking $M=\free{R}i_{*}Y$ in the first spectral
sequence of Theorem~\ref{uct} gives the first spectral sequence
in Theorem~\ref{kt} below.  There is a similar canonical isomorphism
\[
\derF_{R}(\free{R}i_{*}X,F(i_{*}Y,R)) \iso
F(i_{*}X,F(i_{*}Y,R))\iso
F(i_{*}(X\sma Y),R),
\]
in the equivariant stable category.  Via this isomorphism, the second spectral 
sequence of Theorem~\ref{uct} for $M=F(i_{*}Y,R)$ gives the second spectral sequence below.

\begin{thm}(K\"unneth Theorem)
\label{kt}
Let $X$ and $Y$ be $G$-spectra indexed on the same universe.
There is a natural strongly convergent homology spectral sequence of 
$\pR$-modules
\[
E^{2}_{s,\tindex}=\MTor^{\pR}_{s,\tindex}(\hmac{R}_{*} X,\hmac{R}_{*} Y) 
\quad \Longrightarrow \quad 
\hmac{R}_{s+\tindex}(X\sma Y),
\]
and a natural conditionally convergent cohomology spectral sequence of $\pR$-modules
\[
E_{2}^{s,\tindex}=\MExt_{\pcR}^{s,\tindex}(\hmac{R}_{-*} X, \hmac{R}^{*}Y) 
\quad \Longrightarrow \quad 
\hmac{R}^{s+\tindex}(X\sma Y).
\]
\end{thm}

Plugging $M=\free{R}X$ into the spectral
sequences of Theorem~\ref{auct} yields the following spectral sequences.

\begin{thm}\label{akt}
Let $X$ be a $G$-spectrum and let $E$ be a $G$-spectrum indexed on the
complete universe.
There is a natural strongly convergent homology spectral sequence of
$\pR$-modules
\[
E^{2}_{s,\tindex}=\MTor^{\pR}_{s,\tindex}(\hmac{E}_{*} R,\hmac{R}_{*} X) 
\quad \Longrightarrow \quad 
\hmac{E}_{s+\tindex}(R\sma X),
\]
and a natural conditionally convergent cohomology spectral sequence of
$\pR$-modules 
\[
E_{2}^{s,t}=\MExt_{\pcR}^{s,t}(\hmac{R}_{-*} X, \hmac{E}^{*}R) 
\quad \Longrightarrow \quad 
\hmac{E}^{s+t}(R\sma X).
\]
\end{thm}

The spectral sequences above can be combined with Spanier-Whitehead
duality to obtain additional spectral sequences.  If $X$ is a finite
$G$-spectrum, then the Spanier-Whitehead dual of $i_{*}X$ is a
$G$-spectrum $DX$ (indexed on the complete universe).  There are canonical
isomorphisms $E_{*}DX\iso E^{-*}X$ and $E^{*}DX\iso E_{-*}X$, natural in
$E$ and in $X$.  Plugging $DX$ into the spectral sequences of
Theorem~\ref{uct} gives us cohomology to cohomology and cohomology 
to homology universal coefficient spectral sequences.

\begin{thm}
(Dual Universal Coefficient Spectral Sequence)
Let $X$ be a finite $G$-spectrum, and $M$ a
left $R$-module. 
There is a natural strongly convergent homology spectral sequence of 
$\pR$-modules
\[
E^{2}_{s,\tindex}=\MTor^{\pcR}_{s,\tindex}(\hmac{R}^{*} X,\pcM)
\quad \Longrightarrow \quad 
\hmac{M}^{-(s +\tindex)}X
\]
and a natural conditionally convergent cohomology spectral sequence of 
$\pR$-modules
\[
E_{2}^{s,\tindex}=\MExt_{\pR}^{s,\tindex}(\hmac{R}^{-*}X,\pM) 
\quad \Longrightarrow \quad 
\hmac{M}_{-(s +\tindex)}X.
\]
\end{thm}

These dual universal coefficient spectral sequences imply the
following K\"unneth spectral sequences. 

\begin{thm}(Dual K\"unneth Theorem)
Let $X$ and $Y$ be $G$-spectra indexed on the same universe, and
assume that $X$ is finite.
There is a natural strongly convergent homology spectral sequence of 
$\pR$-modules
\[
E^{2}_{s,\tindex}=\MTor^{\pcR}_{s,\tindex}(\hmac{R}^{*} X,\hmac{R}^{*} Y) 
\quad \Longrightarrow \quad 
\hmac{R}^{-(s +\tindex)}(X\sma Y),
\]
and a natural conditionally convergent cohomology spectral sequence of 
$\pR$-modules
\[
E_{2}^{s,\tindex}=\MExt_{\pR}^{s,\tindex}(\hmac{R}^{-*} X, \hmac{R}_{*}Y) 
\quad \Longrightarrow \quad 
\hmac{R}_{-(s +\tindex)}(X\sma Y).
\]
\end{thm}


\section{Graded Mackey functors}\label{secgrdmackey}

This section introduces the categories of \ROG-graded and \Z-graded
Mackey functors. Box products and function objects for
these categories are defined in terms of box products and function
objects for ungraded Mackey functors.  In this discussion, some
familiarity with (ungraded) Mackey functors is assumed.  However,
certain key definitions are reviewed in order to fix notation or to
explain our perspective.

Recall that the Burnside category \Burn may be defined as the
essentially small additive category whose objects are the
finite $G$-sets and whose abelian group of morphisms between the
finite $G$-sets $X$ and $Y$ is the abelian group 
\[
\Burn(X,Y)
=
\shom{\Sigma^\infty
X_{+}}{\Sigma^\infty Y_{+}},
\]
of morphisms between the associated suspension
spectra in the equivariant stable category.  Disjoint union
of finite $G$-sets provides the direct sum operation giving \Burn its
additive structure.
A rather simple, purely algebraic description of
the morphisms in \Burn can be found in \cite[\S V.9]{V1213};
however, 
the approach to the category of Mackey
functors best suited to our purposes is in 
terms of the stable homotopy theoretic description of \Burn.   

\begin{defn}
The category \Mack of Mackey functors is the category of contravariant
additive functors from the Burnside category \Burn to the category \AB
of abelian groups.
\end{defn}

Spanier--Whitehead duality provides an isomorphism of categories
between \Burn and $\Burn^{\op}$ that is the identity on the objects.  
The category of Mackey functors could therefore be defined
equivalently as the category of covariant functors from \Burn to \AB.
Mackey functors can also be described in terms of orbits.  Since
every finite $G$-set is a disjoint union of orbits, every object of
\Burn is isomorphic to a direct sum of the canonical orbits $G/H$.  Let 
\OBurn be the full subcategory of \Burn whose objects are the orbits 
$G/H$ associated to the subgroups $H$ of $G$.  The inclusion of \OBurn into
\Burn induces an equivalence between the category of contravariant
additive functors from \OBurn into \AB and the category \Mack of
Mackey functors.  Thus, the category of Mackey functors could be 
defined equivalently as the category of contravariant additive
functors from \OBurn to \AB. Other equivalent definitions may be found 
in \cite[\S V.9]{V1213} and \cite{JTh90}. 

\begin{defn}
An \ROG-graded Mackey functor \mgMM consists of a
Mackey functor \mgMM[\tindex] for each $\tindex\in \ROG$.  A map of
\ROG-graded Mackey functors $\mgMM\to \mgMN$ 
consists of a map of Mackey functors $\mgMM[\tindex]\to
\mgMN[\tindex]$ for each $\tindex\in \ROG$.  \Z-graded Mackey
functors and maps of \Z-graded
Mackey functors are defined analogously.  The categories of
\ROG-graded and \Z-graded Mackey functors are denoted by \RMack and \ZMack,
respectively.
\end{defn}

The ``shift'' functor $\Sigma^{\aindex}$ takes a graded Mackey 
functor \mgMM to the graded Mackey functor given by 
\begin{equation*}
\left(\Sigma^\aindex \mgMM \right)_\tindex = \mgMM[\tindex - \aindex],
\end{equation*}
for $\tindex \in \ROG$ (or $\tindex\in\Z$).

The categories of \ROG-graded Mackey functors and \Z-graded Mackey
functors have all limits and colimits; these are formed
degree-wise and object-wise.  In other words, for any diagram 
$d\mapsto \mgMM(d)$ of graded Mackey functors,
\[
\left(\Colim_{d}\mgMM(d)\right)_{\tindex}(X) = \Colim_{d}
\left(\mgMM[\tindex](d)(X)\right),
\]
and likewise for the limit.  

\begin{prop}
The categories of \ROG-graded and \Z-graded Mackey functors are
complete and cocomplete abelian categories satisfying $AB5$.
\end{prop}

The box product of graded Mackey functors is constructed from the box
product of Mackey functors, which we now review briefly.
The cartesian product of finite $G$-sets together with the canonical
isomorphism 
$\Sigma^\infty X_{+} \sma \Sigma^\infty Y_{+} \iso \Sigma^\infty(X\times Y)_{+}$ 
in the equivariant stable category provides $\Burn$ with a symmetric
monoidal structure.  The unit for this structure is the one-point $G$-set 
\ePt.  The category \Mack inherits a
symmetric monoidal closed structure from 
\Burn.  This product is most easily
described by the coend formula
\begin{equation*}
(\mMM \boxprod \mMN)(X) = \int^{Y,Z \in \Burn} \mMM(Y) \otimes \mMN(Z)
\otimes \Burn(X, Y \times Z).  
\end{equation*}
This box product is characterized by the universal property that maps of Mackey functors from
$\mMM\boxprod \mMN$ to $\mMP$ are in one-to-one correspondence with
natural transformations
of contravariant functors $\Burn\times
\Burn\to \AB$
from $\mMM\otimes \mMN$ to $\mMP\circ \times$.
The unit for the box product is the Burnside ring Mackey functor
$\mMBr=\Burn(-,\ePt)$.  

The internal $\Hom$ functor adjoint to the box product is
denoted \func{\mMM}{\mMN} and is given by the end formula
\begin{equation*}
\func{\mMM}{\mMN}(X) = \int_{Y,Z \in \Burn} \Hom(\mMM(Y) \otimes
\Burn(Z,X \times Y), \mMN(Z)).
\end{equation*}
It may also be computed using the formulae
\begin{equation*}
\func{\mMM}{\mMN}(X) \iso \Mack(\mMM,\mMN_X) \iso \Mack(\mMM_X,\mMN),
\end{equation*}
where $\mMN_{X}$ denotes the Mackey functor $\mMN(-\times X)$.
Day's general results \cite{CCatF} about monoidal structures on
functor categories specialize to give that the box product $\boxprod$ 
and the internal function object
\func{-}{-} provide the category \Mack with a closed symmetric monoidal
structure.

\begin{defn}
Let \mgMM and \mgMN be \ROG-graded Mackey functors.
Define the
\ROG-graded Mackey functor $\mgMM \gbxprd \mgMN$ by
\begin{equation*}
\left(\mgMM \gbxprd \mgMN \right)_\tindex = \bigoplus_{\aindex +
\bindex = \tindex} \mgMM[\aindex] \boxprod \mgMN[\bindex].
\end{equation*}
Define the \ROG-graded Mackey functor \ghom{\mgMM}{\mgMN} by
\begin{equation*}
\ghom[\tindex]{\mgMM}{\mgMN} = \prod_{\bindex - \aindex = \tindex}
\func{\mgMM[\aindex]}{\mgMN[\bindex]}.
\end{equation*}
For \Z-graded Mackey functors, the functors $(-)\zgbxprd(-)$ and
$\ghom[\#]{-}{-}$ are defined analogously.
\end{defn}

Graded box products and graded function objects provide the
categories \RMack and \ZMack with symmetric monoidal closed
structures. The unit is the \ROG-graded (or \Z-graded) Mackey functor
\mgMBr which is \mMBr in degree zero and zero in all other degrees.  
The unit and associativity isomorphisms for the graded box products
are induced by the unit and associativity isomorphisms for the
ordinary box product.  Further, the adjunction relating graded box
products and graded function objects is easily obtained from the
adjunction relating ordinary box products and function objects.  This
describes the monoidal closed structure on \RMack and \ZMack.  For
\Z-graded Mackey functors, the symmetry isomorphism is just as
straightforward.  On the summands of $(\mzgMM
\zgbxprd \mzgMN)_t$ and $(\mzgMN \zgbxprd \mzgMM)_t$ 
for $t = m+n$, this isomorphism is the composite of the Mackey functor
symmetry isomorphism $\mzgMM[m] \boxprod \mzgMN[n] \iso \mzgMN[n]
\boxprod \mzgMM[m]$ with multiplication by $(-1)^{mn}$.

For \ROG-graded Mackey functors, the symmetry isomorphism is more
delicate because more complicated ``signs'' are needed for
compatibility with our topological applications.  These signs are units in the 
Burnside ring $\mMBr(G/G)$ of $G$.  Any element $a$ in the Burnside
ring of $G$ defines a natural transformation in the category of Mackey
functors from the identity functor to itself.  The Yoneda lemma
identifies $\mMBr(G/G)$ as the endomorphism ring of \mMBr, and the
natural transformation on a Mackey functor \mMM associated to $a$ is
the composite
\begin{equation*}
\mMM \iso \mMM \boxprod \mMBr \lrarrow[\id_{\mMM}\boxprod a] \mMM
\boxprod \mMBr \iso \mMM.
\end{equation*}
When $a$ is a unit in $\mMBr(G/G)$, this natural transformation is an
isomorphism. Each unit $a$ in $\mMBr(G/G)$ satisfies $a^2 =1$ so these
units can be thought of as generalized signs.

Given a function $\sigma$ from $\ROG \times \ROG$ to the group of
units of $\mMBr(G/G)$, we can define a natural isomorphism
\[
c_{\sigma}\colon \mgMM \gbxprd \mgMN \iso \mgMN \gbxprd \mgMM
\]
by taking the map on the summands associated to
$\aindex+\bindex=\tindex$ to be the composite of the symmetry
isomorphism $\mgMM[\aindex]\boxprod \mgMN[\bindex] \iso \mgMN[\bindex]
\boxprod \mgMM[\aindex]$ and the automorphism $\sigma(\aindex,
\bindex)$.  An easy diagram chase shows that $c_{\sigma}$ makes the
box product into a symmetric monoidal product exactly when $\sigma$ is
(anti)symmetric, that is,
\begin{equation*}
\sigma(\aindex,\bindex)\sigma(\bindex,\aindex)=1
\end{equation*}
and bilinear, that is,
\begin{equation*}
\sigma(\aindex + \bindex,\cindex) = \sigma(\aindex,\cindex)
\sigma(\bindex,\cindex).
\end{equation*}

The appropriate choice of $\sigma$ is dictated by the definition of
the homotopy Mackey functor \pimac and depends on the
topological choices in that definition.  Because of the required bilinearity
of $\sigma$, it suffices to specify $\sigma(\aindex,\bindex)$ in only
those cases where \aindex and \bindex are irreducible real
$G$-representations.  In Appendix~\ref{secpistar}, we show that, if 
\aindex and \bindex are non-isomorphic, then we can assume that
$\sigma(\aindex,\bindex) = 1$.  The element $\sigma(\aindex,\aindex)$
of $\mMBr(G/G)$ is the automorphism of $S=\Sigma^{\infty}G/G_{+}$
obtained by stabilizing the multiplication by $-1$ map on $S^{\aindex}$.  
This definition of
$\sigma$ completes the construction of the symmetry isomorphism in the
\ROG-graded case and so finishes the construction of the symmetric
monoidal product.

\begin{prop}
The categories \RMack and \ZMack are closed symmetric monoidal abelian
categories.
\end{prop}


\section{Graded Mackey functor rings and modules}\label{secgrdmods}

This section is devoted to a discussion of rings and modules in the categories of graded
Mackey functors.  In particular, the ``box over \mgMR''
(\rbox) and the ``\mgMR-function'' (\rfunc{\mgMR}{-}{-}) constructions for
modules over a graded Mackey functor ring \mgMR are defined, and their basic 
properties are described.  For simplicity, all definitions are given in the 
context of the category \RMack of \ROG-graded Mackey functors.  However, with 
the obvious notational modifications, these definitions apply equally well 
to the category \ZMack of \Z-graded Mackey functors.  We begin with the (usual) 
definitions of rings and modules in a symmetric monoidal abelian category.

\begin{defn}
An \ROG-graded Mackey functor ring consists of an \ROG-graded
Mackey functor \mgMR together with unit $i\colon \mgMBr \to \mgMR$ and
multiplication $\mu \colon \mgMR \gbxprd \mgMR \to \mgMR$ maps for
which the unit and associativity diagrams
\begin{equation*}
\xymatrix@C=2em{ {\mgMBr \gbxprd \mgMR \ar[r]^-{i \gbxprd \id}
\ar[dr]_{\iso} } & {\mgMR \gbxprd \mgMR \ar[d]^{\mu}}
	& {\mgMR \gbxprd \mgMBr \ar[l]_-{\id \gbxprd i} \ar[dl]^{\iso} } \\
{} & {\mgMR} & {} } \quad \xymatrix@C=2em{ {\mgMR \gbxprd \mgMR
\gbxprd \mgMR \ar[r]^-{\mu \gbxprd \id} \ar[d]_{\id \gbxprd \mu} }
	& {\mgMR \gbxprd \mgMR \ar[d]^-{\mu}}  \\
{\mgMR \gbxprd \mgMR \ar[r]_-{\mu}} & {\mgMR} }
\end{equation*}
commute.  The ring \mgMR is said to be commutative if the symmetry
diagram
\begin{equation*}
\xymatrix{
{\mgMR \gbxprd \mgMR \ar[rr]^-{\iso} \ar[dr]_{\mu} } & {} & {\mgMR \gbxprd \mgMR  \ar[dl]^{\mu} } \\
{} & {\mgMR} & {} }
\end{equation*}
also commutes.  In these diagrams, the unlabeled isomorphisms are the unit
and symmetry isomorphisms of the symmetric monoidal category \RMack.
\end{defn}

\begin{defn}
A left module over an \ROG-graded Mackey functor ring \mgMR
consists of an \ROG-graded Mackey functor \mgMM and an action map $\xi
\colon \mgMR \boxprod \mgMM \to \mgMM$ for which the unit and
associativity diagrams
\begin{equation*}
\xymatrix@C=2em{
{\mgMBr \gbxprd \mgMM \ar[r]^-{i \gbxprd \id} \ar[dr]_{\iso} } & {\mgMR \gbxprd \mgMM \ar[d]^{\xi}}  \\
{} & {\mgMM} } \quad \xymatrix@C=2em{
{\mgMR \gbxprd \mgMR \gbxprd \mgMM \ar[r]^-{\mu \gbxprd \id} \ar[d]_{\id \gbxprd \xi} } & {\mgMR \gbxprd \mgMM \ar[d]^-{\xi}}  \\
{\mgMR \gbxprd \mgMM \ar[r]_-{\xi}} & {\mgMM} }
\end{equation*}
commute.  A right module \mgMM over \mgMR is defined analogously in
terms of an action map $\zeta \colon \mgMM \boxprod \mgMR \to \mgMM$.
The categories of left and right modules over an \ROG-graded Mackey
functor ring \mgMR are denoted \rMod and \Modr, respectively.
\end{defn}

For our study of the additional structure carried by some of our spectral 
sequences, we also need the notion of \mgMR-bimodule.

\begin{defn}
An \mgMR-bimodule is a graded Mackey functor \mgMM having left
and right \mgMR-actions ($\xi$ and $\zeta$, respectively) for which
the diagram 
\begin{equation*}
\xymatrix@C=2em{ {\mgMR \gbxprd \mgMM \gbxprd \mgMR \ar[r]^-{\xi
\gbxprd \id} \ar[d]_{\id \gbxprd \zeta} }
	& {\mgMM \gbxprd \mgMR \ar[d]^-{\zeta}}  \\
{\mgMR \gbxprd \mgMM \ar[r]_-{\xi}} & {\mgMM} }
\end{equation*}
commutes.  If \mgMR is commutative, then every left module \mgMM
carries a bimodule structure in which the right action map $\zeta$ is
the composite
\[
\mgMM \gbxprd \mgMR \iso \mgMR \gbxprd \mgMM \lrarrow[\xi] \mgMM.
\]
A bimodule structure of this form is called symmetric.
\end{defn}

When \mgMR is an \ROG-graded Mackey functor ring and \mgMK is an
\ROG-graded Mackey functor, the \ROG-graded Mackey functors
$\mgMR \gbxprd \mgMK$ and \ghom{\mgMR}{\mgMK} carry \mgMR-bimodule
structures coming from the left and right actions of \mgMR on itself.
Regarded simply as left \mgMR-modules, $\mgMR \gbxprd \mgMK$ and
\ghom{\mgMR}{\mgMK} are called the free and cofree left \mgMR-modules
generated by \mgMK.  The free and cofree right \mgMR-modules generated
by \mgMK are constructed similarly.  A purely formal argument indicates
that free and cofree \mgMR-modules 
have the expected universal properties.

\begin{prop}
Let \mgMR be an \ROG-graded Mackey functor ring.  The functor
$\mgMR\gbxprd(-)\colon 
\RMack \to \rMod$ is left adjoint to the forgetful functor $\rMod \to
\RMack$. The functor $\ghom{\mgMR}{-}\colon \RMack \to \rMod$ is right 
adjoint to the forgetful functor $\rMod \to \RMack$. 
\end{prop}

These two adjunctions may be used to identify the category of left
\mgMR-modules as both the category of algebras over a monad on \RMack
and the category of coalgebras over a comonad on \RMack. These
identifications together with the corresponding identifications for
the category of right \mgMR-modules and the category of
\mgMR-bimodules are the key to proving the following result.

\begin{prop}\label{rModcomplete}
Let \mgMR be an \ROG-graded Mackey functor ring.  The
categories of left \mgMR-modules, right \mgMR-modules, and
\mgMR-bimodules are complete 
and cocomplete abelian categories that satisfy $AB5$.
\end{prop}

We now move on to the main constructions of this section, the functors
$(-)\rbox(-)$ and \rfunc{\mgMR}{-}{-}.

\begin{defn}
Let $\mgML$ and $\mgMM$ be left $\mgMR$-modules and let $\mgMN$ be a
right $\mgMR$-module. 
\begin{thmlist}
\item
The graded Mackey functor $\mgMN\rbox\mgMM$ is defined by the
coequalizer diagram 
\begin{equation*}
\xymatrix{ {\mgMN \gbxprd \mgMR \gbxprd \mgMM \ar@<0.7ex>[r]^-{\zeta
\gbxprd \id } \ar@<-0.7ex>[r]_-{\id \gbxprd \xi}} & {\mgMN \gbxprd
\mgMM \ar[r]^-{q}} & {\mgMN \rbox \mgMM}. }
\end{equation*}
Here $\xi$ is the action map for $\mgMM$ and $\zeta$ is the action map
for $\mgMN$.
\item 
The \ROG-graded Mackey functor \rfunc{\mgMR}{\mgML}{\mgMM} 
is defined by the equalizer diagram
\begin{equation*}
\xymatrix{ {\rfunc{\mgMR}{\mgML}{\mgMM} \ar[r]^-{j}} &
{\ghom{\mgML}{\mgMM} \ar@<0.7ex>[r]^-{{\nu^*}}
\ar@<-0.7ex>[r]_-{\tilde\xi_*}} & {\ghom{\mgMR \gbxprd \mgML}{\mgMM} \iso
\ghom{\mgML}{\ghom{\mgMR}{\mgMM}}.}  }
\end{equation*}
Here the map $\nu^*$ is induced from the action map $\nu$ for \mgML,
and the map $\tilde\xi_*$ is induced from the coaction map
$\tilde\xi\colon{\mgMM}\to{\ghom{\mgMR}{\mgMM}}$ adjoint to the action
map for \mgMM.
\end{thmlist}
\end{defn}

The following proposition follows directly from these definitions and the
properties of $\gbxprd$ and \ghom{-}{-}.

\begin{prop}\label{RModHalfExact} 
Let $\mgMM$ and $\mgMN$ be a left and a right $\mgMR$-module, respectively. 
\begin{thmlist}
\item The functors 
${\mgMN \rbox (-)}\colon{\rMod}\to{\RMack}$ and ${(-) \rbox
\mgMM}\colon{\Modr}\to{\RMack}$
preserve all colimits and are therefore right exact.

\item The functor 
${\rfunc{\mgMR}{\mgMM}{-}}\colon{\rMod}\to{\RMack}$
preserves all limits and is therefore left exact. 

\item The functor 
${\rfunc{\mgMR}{-}{\mgMM}}\colon{\rMod^{\op}}\to{\RMack}$
converts colimits in \rMod into limits in \RMack and is therefore left exact.
\end{thmlist}
\end{prop} 

For any graded Mackey functors \mgMM, \mgMMp, and \mgMMpp, 
the closed symmetric monoidal structure on \RMack provides a
composition pairing  
\begin{equation*}
\ghom{\mgMMp}{\mgMMpp} \gbxprd \ghom{\mgMM}{\mgMMp} \lrarrow
\ghom{\mgMM}{\mgMMpp} 
\end{equation*}
that is both associative and unital.  This pairing restricts to 
a pairing for the category of left \mgMR-modules.  

\begin{prop}\label{CompPair}
Let \mgMM, \mgMMp, and \mgMMpp be left \mgMR-modules.  The 
composition pairing on \ghom{-}{-} restricts to a pairing 
\begin{equation*}
\rfunc{\mgMR}{\mgMMp}{\mgMMpp} \gbxprd \rfunc{\mgMR}{\mgMM}{\mgMMp}
\lrarrow \rfunc{\mgMR}{\mgMM}{\mgMMpp}
\end{equation*}
that is both associative and unital.
\end{prop}

Although the constructions $\mgML \rbox \mgMM$ and
\rfunc{\mgMR}{\mgMM}{\mgMN} typically yield only \ROG-graded Mackey
functors, in some important special cases they
yield \mgMR-modules. 

\begin{prop}
Let \mgML and \mgMM be left \mgMR-modules, and let \mgMN be a right
\mgMR-module. 
\begin{thmlist}
\item If \mgML is an \mgMR-bimodule, then 
\rfunc{\mgMR}{\mgML}{\mgMM} is naturally a left \mgMR-module.

\item If \mgMM is an \mgMR-bimodule, then $\mgMN \rbox \mgMM$ and
\rfunc{\mgMR}{\mgML}{\mgMM} are naturally right \mgMR-modules.

\item If \mgMN is an \mgMR-bimodule, then $\mgMN \rbox \mgMM$ is
naturally a left \mgMR-module.
\end{thmlist}
\end{prop}

When \mgMR is commutative, we can always consider the symmetric
\mgMR-bimodule structure on any (left or right) \mgMR-module to obtain
\mgMR-module structures on $(-)\rbox(-)$ and \rfunc{\mgMR}{-}{-}.  In
fact, regarding either \mgMM or \mgMN as the bimodule yields the same
\mgMR-module structure on $\mgMM \rbox \mgMN$, and similarly for
\rfunc{\mgMR}{\mgMM}{\mgMN}.   In this context the following stronger
version of the previous results hold.

\begin{prop}\label{RComm}
Let \mgMR be a commutative graded Mackey functor ring.
\begin{thmlist}
\item
The module category \rMod is a closed symmetric monoidal abelian
category with product \rbox and function object
\rfunc{\mgMR}{-}{-}.
\item
The free functor
${\mgMR\gbxprd(-)}\colon {\RMack}\to{\rMod}$ is strong symmetric monoidal
and the forgetful functor $\rMod\to\RMack$ is lax symmetric
monoidal. 
\item
The composition pairing 
\[
\rfunc{\mgMR}{\mgMMp}{\mgMMpp} \rbox \rfunc{\mgMR}{\mgMM}{\mgMMp}
\lrarrow \rfunc{\mgMR}{\mgMM}{\mgMMpp} 
\]
coming from the closed symmetric monoidal structure on \rMod is the obvious one derived from the composition pairing of Proposition~\ref{CompPair}.  

\end{thmlist}
\end{prop}


\section{The homological algebra of graded Mackey functor modules}\label{sechom}

This section is devoted to homological algebra for the categories of
Mackey functor modules over a graded Mackey functor ring.  Our first 
objective is to show that these categories have enough projectives 
and injectives. These objects are then used to construct the 
derived functors $\MTor^{\mgMR}_*$ and $\MExt^*_{\mgMR}$ in terms of
resolutions and to show that they have the expected properties.  Some 
notation is needed to construct the desired injective and projective objects.  

\begin{defn}
For each finite $G$-set $X$, let $\mMBr^{X}$ denote the
Mackey functor $\Burn(-,X)$.  For any abelian group $E$, let  
$\inj{X}{E}$ denote the Mackey functor $\Hom(\Burn(X,-),E)$.  The 
corresponding graded Mackey functors concentrated in degree zero are 
denoted $\mgMBr^{X}$ and $\ginj{X}{E}$.
\end{defn}

The enriched Yoneda Lemma gives natural isomorphisms 
\begin{alignat*}{2}
\Mack(\mMBr^{X},\mMM) &\iso \mMM(X) & \qquad 
\Mack(\mMM, \inj{X}{E}) &\iso \Hom(\mMM(X),E) \\
\func{\mMBr^{X}}{\mMM}(Y) & \iso \mMM(X\times Y) & \qquad
\func{\mMM}{\inj{X}{E}}(Y) &\iso \Hom(\mMM(X\times Y),E)
\end{alignat*}
of abelian groups.  A coend argument dual to the end argument 
used to prove the Yoneda Lemma provides the well-known isomorphism
\[
(\mMBr^{X} \boxprod \mMM)(Y)\iso \mMM(X\times Y).
\]
In the graded context, these isomorphisms give the following proposition:  

\begin{prop}\label{grflat}
Let $\mgMM$ be a graded Mackey functor.
There are natural isomorphisms of abelian groups
\begin{gather*}
\ghom[\aindex]{\Sigma^{\tindex}\mgMBr^{X}}{\mgMM}(Y)\iso 
\mgMM[\tindex +\aindex](X\times Y) \\
(\Sigma^{\tindex}\mgMBr^{X}\gbxprd \mgMM)_{\aindex}(Y)\iso
\mgMM[-\tindex+\aindex](X\times Y)\\
\ghom[\aindex]{\mgMM}{\Sigma^{\tindex} \ginj{X}{E}}(Y) \iso 
\Hom(\mgMM[\tindex -\aindex](X\times Y),E).
\end{gather*}
In particular, $\Sigma^{\tindex}\mgMBr^{X}$ is a projective object in
\RMack and 
the functors $\ghom{\Sigma^{\tindex}\mgMBr^{X}}{-}$ 
and $\mgMBr^{X}\gbxprd (-)$ are exact.  Also, if $E$ is an injective abelian
group, then $\Sigma^{\tindex}\ginj{X}{E}$ is an injective object in \RMack, and the
functor $\ghom[]{-}{\Sigma^{\tindex}\ginj{X}{E}}$ is 
exact. 
\end{prop}

It follows from this proposition that the category of
graded Mackey functors has enough projectives and injectives.  
Let \mgMM be a graded Mackey functor. For each $\tindex$ in
\ROG and each subgroup $H$ of $G$, choose a surjection
$P_{\tindex,H}\to \mgMM[\tindex](G/H)$ whose domain is a free abelian 
group and an injection $\mgMM[\tindex](G/H)\to I_{\tindex,H}$ whose 
range is an injective abelian group $I_{\tindex,H}$.  Then the
previous proposition provides an epimorphism
\[
\mgMP=\bigoplus_{\tindex ,H}\Sigma^{\tindex}\mgMBr^{G/H}\otimes
{P_{\tindex,H}}\lrarrow \mgMM 
\]
and a monomorphism
\[
\mgMM \lrarrow 
\prod_{\tindex,H}\Sigma^{\tindex}\ginj{G/H}{I_{\tindex,H}}=\mgMI.
\]
When $\mgMM$ is a left $\mgMR$-module for some graded Mackey functor ring
$\mgMR$, then the induced epimorphism 
$\mgMR\gbxprd \mgMP\to \mgMM$ of left $\mgMR$-modules has domain 
a projective left
$\mgMR$-module.  Likewise, the induced monomorphism 
$\mgMM\to \ghom{\mgMR}{\mgMI}$ of left $\mgMR$-modules has codomain 
an injective left
$\mgMR$-module.  Similar observations apply in the case of right
modules.  Thus, we have proven:  

\begin{prop}
Let $\mgMR$ be a graded Mackey functor ring.  The categories of left and
right $\mgMR$-modules have enough projectives and injectives.
\end{prop}

Since an epimorphism from one projective onto another or a
monomorphism from one injective into another is split, the preceding
argument also provides characterizations of projective and injection modules.

\begin{prop}\label{charproj}
An $\mgMR$-module is projective if and only if it is a direct
summand of a direct sum of $\mgMR$-modules of the form $\mgMR\gbxprd 
\Sigma^{\tindex}\mgMBr^{G/H}$.  Also, an $\mgMR$-module is injective 
if and only if it is
a direct summand of a product of $\mgMR$-modules of the form
$\ghom{\mgMR}{\Sigma^{\tindex}{\ginj{G/H}{I}}}$ in which $I$ is an
injective abelian group.
\end{prop}

Proposition~\ref{grflat} and this characterization of
projectives and injectives have some important implications 
for the exactness of the functors $\rbox$ and $\rfunc{\mgMR}{-}{-}$.  
In the terminology of Lewis~\cite{PrjFlt},
these exactness results assert that projectives and injectives are respectively
internal projective and internal injective and that projective implies
flat.

\begin{thm}\label{homological}
Let $\mgMR$ be a Mackey functor ring.
\begin{thmlist}
\item If $\mgMP$ is a projective left $\mgMR$-module, then
$\rfunc{\mgMR}{\mgMP}{-}$ is an exact functor from left $\mgMR$-modules
to graded Mackey functors.
\item If $\mgMI$ is an injective left $\mgMR$-module, then 
$\rfunc{\mgMR}{-}{\mgMI}$ is an exact functor from left $\mgMR$-modules
to graded Mackey functors.
\item If $\mgMP$ is a projective left $\mgMR$-module, then
$(-)\rbox\mgMP$ is an exact functor from right $\mgMR$-modules to
graded Mackey functors.
\item If $\mgMQ$ is a projective right $\mgMR$-module, then
$\mgMQ\rbox(-)$ is an exact functor from left $\mgMR$-modules to
graded Mackey functors.
\end{thmlist}
\end{thm} 

In the category of left \mgMR-modules, as in any abelian category 
with enough projectives and injectives, every object has both projective and
injective resolutions.  To be consistent with
the usual grading conventions for spectral sequences, we use the inner
degree for the resolution degree in projective resolutions.  Thus,
for any projective resolution \mbMP{*}{*}, \mbMP{s}{*} is an
\mgMR-module for each $s\geq 0$ and
\mbMP{*}{\tindex} is a chain complex of Mackey functors for each
$\tindex\in \ROG$ (or \Z).  When \mcMI{*}{*} is an injective
resolution, \mcMI{s}{*} is an \mgMR-module for each $s\geq 0$ and
\mcMI{*}{\tindex} is a cochain complex of 
Mackey functors for each $\tindex\in \ROG$ (or \Z).

If \mbMP{*}{*} is a projective resolution of a left \mgMR-module
\mgMM and \mbMQ{*}{*} is a projective resolution of a right
\mgMR-module \mgMN, then the maps
\begin{equation*}
\mgMN \rbox \mbMP{*}{*} \llarrow \mbMQ{*}{*} \rbox \mbMP{*}{*}
\lrarrow \mbMQ{*}{*} \rbox \mgMM
\end{equation*}
are quasi-isomorphisms of chain complexes of graded Mackey functors by
Theorem~\ref{homological}. The graded Mackey functor
\mTor[s]{\mgMR}{\mgMN}{\mgMM} is defined to be the $s$-th homology of any of
these cochain complexes.  The Mackey functor 
$(\mTor[s]{\mgMR}{\mgMN}{\mgMM})_{\tindex}$ is usually denoted $\mTor[{s,\tindex}]{\mgMR}{\mgMN}{\mgMM}$.  

If \mbMO{*}{*} is a projective resolution of a left \mgMR-module
\mgML and \mcMI{*}{*} is an injective resolution of a
left \mgMR-module \mgMM, then the maps
\begin{equation*}
\rfunc{\mgMR}{\mbMO{*}{*}}{\mgMM} \lrarrow
\rfunc{\mgMR}{\mbMO{*}{*}}{\mcMI{*}{*}} \llarrow
\rfunc{\mgMR}{\mgML}{\mcMI{*}{*}}
\end{equation*}
are quasi-isomorphisms of cochain complexes of graded Mackey functors
by Theorem~\ref{homological}. The graded Mackey functor
\mExt[s]{\mgMR}{\mgML}{\mgMM} is defined to be the $s$-th cohomology of any of
these cochain complexes.  The Mackey functor 
$(\mExt[s]{\mgMR}{\mgML}{\mgMM})_{-\tindex}$ is usually denoted $\mExt[{s,\tindex}]{\mgMR}{\mgML}{\mgMM}$.  

The standard homological arguments imply that $\MTor^{\mgMR}_*$ 
and $\MExt^*_{\mgMR}$ behave as one would expect.

\begin{prop}\label{TorExtProp}
The constructions $\MTor^{\mgMR}_*$ and $\MExt^*_{\mgMR}$ are
well-defined, natural in each variable, and convert short exact
sequences in each variable to long exact sequences.  Moreover, there
are canonical natural isomorphisms
\begin{equation*}
\mTor[0]{\mgMR}{\mgML}{\mgMM} \iso \mgML \rbox \mgMM
\end{equation*}
and
\begin{equation*}
\mExt[0]{\mgMR}{\mgMM}{\mgMN} \iso \rfunc{\mgMR}{\mgMM}{\mgMN}.
\end{equation*}
\end{prop}

As with {\rbox} and \rfunc{\mgMR}{-}{-},
if one of the arguments of $\MTor^{\mgMR}_*$ or $\MExt^*_{\mgMR}$ is
an \mgMR-bimodule, then $\MTor^{\mgMR}_*$ or $\MExt^*_{\mgMR}$
inherits an \mgMR-module structure from that bimodule.  In the 
context of ordinary rings, this assertion is a
purely formal consequence of the functoriality of $\MTor^{\mgMR}_*$
and $\MExt^*_{\mgMR}$.  However, since \mgMR is a graded Mackey
functor ring, to obtain this result, we must either show that the
functors $\MTor^{\mgMR}_*$ and $\MExt^*_{\mgMR}$ are enriched over
\RMack or provide a direct construction of
the \mgMR-action.  Either of these approaches requires the following
lemma.  Its proof, which uses the right exactness of
$\rbox$, is essentially identical to that of the corresponding result
for chain complexes of abelian groups. 

\begin{lem}\label{HcHpairings}
Let \mbMC{*}{*} and \mbMD{*}{*} be chain complexes of graded Mackey
functors, and \mcMC{*}{*} and \mcMD{*}{*} be cochain complexes of
graded Mackey functors.  Then there are natural transformations
\begin{equation*}
H_*(\mbMC{*}{*}) \gbxprd H_*(\mbMD{*}{*}) \lrarrow H_*(\mbMC{*}{*}
\gbxprd \mbMD{*}{*})
\end{equation*}
and
\begin{equation*}
H^*(\mcMC{*}{*}) \gbxprd H^*(\mcMD{*}{*}) \lrarrow H^*(\mcMC{*}{*}
\gbxprd \mcMD{*}{*})
\end{equation*}
which are unital and associative in the appropriate sense.
\end{lem}

Taking one of the complexes in this lemma to be \mgMR
(concentrated in homological or cohomological degree zero) and the
other to be the complex used to compute $\MTor^{\mgMR}_{*}$ or
$\MExt^{*}_{\mgMR}$, we obtain the desired actions of \mgMR.

\begin{thm}\label{TorBiM}
Let \mgML and \mgMM be left
\mgMR-modules, and let \mgMN be a right \mgMR-module.
\begin{thmlist}
\item If \mgML is an \mgMR-bimodule, then \mExt[s]{\mgMR}{\mgML}{\mgMM} 
is naturally a right \mgMR-module. 

\item If \mgMM is an \mgMR-bimodule, then
\mTor[s]{\mgMR}{\mgMN}{\mgMM} and
\mExt[s]{\mgMR}{\mgML}{\mgMM} are naturally right \mgMR-modules.

\item If \mgMN is an \mgMR-bimodule, then
\mTor[s]{\mgMR}{\mgMN}{\mgMM} is naturally a left \mgMR-module.

\end{thmlist} 
\end{thm}

The second natural transformation in Lemma~\ref{HcHpairings} may also be used to 
construct the Yoneda pairing for $\MExt^*_{\mgMR}$.  Let \mgMM,
\mgMMp, and \mgMMpp be left \mgMR-modules, \mbMP{*}{*} be a projective
resolution of \mgMM, and \mcMI{*}{*} be an injective resolution of
\mgMMpp.  Then the Yoneda pairing
\begin{equation*}
\mExt{\mgMR}{\mgMMp}{\mgMMpp} \gbxprd \mExt{\mgMR}{\mgMM}{\mgMMp}
\lrarrow \mExt{\mgMR}{\mgMM}{\mgMMpp}
\end{equation*}
is the composite
\begin{eqnarray*}
H^{*}(\rfunc{\mgMR}{\mgMMp}{\mcMI{*}{*}}) \boxprod
H^{*}(\rfunc{\mgMR}{\mbMP{*}{*}}{\mgMMp}) 
& \lrarrow &
H^{*}(\rfunc{\mgMR}{\mgMMp}{\mcMI{*}{*}} \gbxprd 
\rfunc{\mgMR}{\mbMP{*}{*}}{\mgMMp}) \\
 & \lrarrow & H^{*}(\rfunc{\mgMR}{\mbMP{*}{*}}{\mcMI{*}{*}}).  
\end{eqnarray*}
Here, the first map is the second natural transformation from 
the lemma and the second map comes from the composition pairing 
of Lemma \ref{CompPair}.  The usual homological arguments imply 
that this pairing on $\MExt^*_{\mgMR}$ behaves as expected.

\begin{prop}\label{Ypair}
The Yoneda pairing is well-defined and associative.  Moreover, in
cohomological degree zero, it agrees with the usual composition
pairing.
\end{prop}


\section{Homotopy Mackey functors and the homological algebra of
equivariant $R$-modules}\label{sectop} 

This section is devoted to a discussion of the functor $\pimac$ 
from the equivariant stable 
category to the category of graded Mackey functors.  This functor 
connects the derived smash product and function object constructions 
for modules over an equivariant $S$-algebra $R$ to the corresponding
$\rbox[\pimac R]$ and $\rfunc{\pimac R}{-}{-}$ constructions.  As a 
first step in our analysis of the link between the 
homotopy theory of $R$-modules and the homological
algebra of $\pimac R$-modules, we consider here those $R$-modules 
whose homotopy Mackey functors are projective or
injective as $\pimac R$-modules.  

We begin by reviewing the construction of $\pimac$.  The first 
step in this process is selecting a model $S^{\tindex}$ of 
the \tindex-sphere for each element \tindex of \ROG.  It is convenient to take
$S^{0}=\gsS$.  More generally, when \tindex is the trivial
representation of dimension $n>0$, $S^{\tindex}$ is taken to be the 
smash product of $\gsS$ with the standard $n$-sphere space. For all
other \tindex, the object $S^{\tindex}$ may be chosen arbitrarily from
the appropriate homotopy class.  Then, \pimac is defined for any
$G$-spectrum $\gsM$ and any finite $G$-set $X$ by
\begin{equation*}
\pimacS[\tindex]{\gsM}(X) = \shom{S^{\tindex}\sma X_{+}}{\gsM}. 
\end{equation*}
The naturality of this construction in stable maps of $X$ gives
\pimacS[\tindex]{\gsM} the structure of a Mackey functor.  Further,
its naturality in \gsM makes \pimac a functor from the equivariant
stable category to \RMack.  Restricting the index \tindex so that it
lies in \Z rather than \ROG gives a functor, which we also denote
\pimac, from the equivariant stable category to \ZMack.  

Appendix~\ref{secpistar} contains a complete proof of the following 
folk theorem:  

\begin{thm}\label{pisym}
The functor \pimac is a lax symmetric
monoidal functor from the equivariant stable category to the category
of \ROG-graded (or \Z-graded) Mackey functors.  
\end{thm}

In other words, we have a suitably associative, symmetric, and unital
natural transformation $\pimac \gsN\gbxprd \pimac \gsM\to
\pimac(\gsN\sma \gsM)$.  By comparing the diagrams used to define 
homotopical ring and module spectra in the equivariant stable category 
with those used to define graded Mackey functor rings and modules, it 
is easy to see that $\pimac$
takes homotopical ring and module spectra to graded Mackey functor rings and
modules.  We apply this to equivariant $S$-algebras and modules over an 
equivariant $S$-algebra $R$ in a modern category of spectra, which pass to 
homotopical ring and module spectra in the equivariant stable category.  
A weak bimodule in the category of left $R$-modules is defined to be left $R$-module 
together with
a homotopical right $R$-module structure in the derived category of
left $R$-modules.  Clearly, the underlying spectrum of a weak bimodule is 
a homotopical bimodule.  Weak bimodules in the category of right
$R$-modules are defined analogously; their  
underlying spectra are also homotopical bimodules.  These observations are summarized in the following 
corollary of Theorem \ref{pisym}:

\begin{cor}
Let $R$ be an equivariant $S$-algebra.  Then $\pimac R$ is a graded
Mackey functor ring which is commutative if $R$ is.  The functor 
$\pimac$ refines
to a functor from the category of left (or right) $R$-modules to the 
category of left (or right) $\pimac
R$-modules.  Both refinements take weak bimodules to 
$\pimac R$-bimodules.
\end{cor}

Let $M$ and $N$ be left and right modules, respectively, over 
an equivariant $S$-algebra $R$.  Recall that the functor $\dersma^{R}$ is the 
derived functor of
$\sma_{R}$ and that $\TTor^{R}_{*}(N,M)$ denotes $\pimac
\dersma^{R}(N,M)$.  There is a canonical comparison map 
\[
N\sma M \to \dersma^{R}(N,M) 
\]
derived from the natural map $N \sma_{S} M \to N \sma_{R} M$. 
The two composites 
\[
N\sma R\sma M\to N \sma  M \to \dersma^{R}(N,M) 
\]
coming from the actions of $R$ on $M$ and $N$ coincide.  Thus, the 
universal property defining $\rbox[\pimac R]$ provides a natural transformation 
\[
\pimac N\rbox[(\pimac R)] \pimac M\lrarrow \TTor^{R}_{*}(N,M).
\]
There is a canonical natural
isomorphism  
\[
\dersma^{R}(N,M)\sma X\lrarrow \dersma^{R}(N,M\sma X)
\]
in the equivariant stable category.  This isomorphism is associative in  
the obvious sense and makes the diagram 
\[
\xymatrix@C=-1em{
&N\sma M\sma X\ar[dr]\ar[dl]\\
\dersma^{R}(N,M)\sma X\ar[rr]&&\dersma^{R}(N,M\sma X)
}
\]
commute.  There is an analogous isomorphism in the other variable with analogous
properties.  These isomorphisms, and their asserted properties, come from the fact that 
$\dersma^{R}(N,-)$ is enriched functor over the equivariant stable category
\cite{LMmmmc}.  Taking $X=R$, we see that, 
when $M$ or $N$ is a weak bimodule, $\dersma^{R}(N,M)$ is naturally a
homotopical $R$-module and the comparison map $N\sma M \to
\dersma^{R}(N,M)$ is a map of homotopical
$R$-modules.  The implications of these observations are summarized 
in the following result:  

\begin{thm}\label{ttormap}
Let $R$ be an equivariant $S$-algebra, $M$ be a left $R$-module, and $N$ 
be a right $R$-module.  There is a natural
transformation of graded Mackey functors
\[
\pimac N\rbox[(\pimac R)] \pimac M\lrarrow \TTor^{R}_{*}(N,M).
\]
If $M$ is a weak bimodule, this is a map of right
$\pimac R$-modules.  If $N$ is an weak bimodule, it is a map of
left $\pimac R$-modules.
\end{thm}

Likewise, the functor $\derF_{R}$ is the derived functor of
$F_{R}$, and $\TExt_{R}^{-*}(L,M)=\pimac \derF_{R}(L,M)$ for 
left $R$-modules $L$ and $M$.    
The natural map $F_{R}(L,M)\to F_{S}(L,M)$ induces a map 
$\derF_{R}(L,M)\to F(L,M)$.  The two maps from 
$\derF_{R}(L,M)$ to $F(R\sma L,M)$ coming from the actions of $R$ on $L$ 
and $M$ coincide and so induce a 
natural transformation 
\[
\TExt_{R}^{-*}(L,M)\lrarrow
\rfunc{(\pimac R)}{\pimac L}{\pimac M}.
\]
The functors $\derF_{R}(L,-)$ and $\derF_{R}(-,M)$ are also enriched
over the equivariant stable 
category, and so we have natural maps
\[
\derF_{R}(L,M)\sma X\to \derF_{R}(L,M\sma X)
\quad\text{and}\quad
\derF_{R}(L\sma X,M)\to F(X,\derF_{R}(L,M))
\]
in the equivariant stable category.  The second of these maps is always 
an isomorphism, but the first map generally
is not.  These maps satisfy the evident associativity conditions and are 
compatible with the canonical maps 
\[
F(L,M)\sma X\to F(L,M\sma X)
\quad\text{and}\quad
F(L\sma X,M)\to F(X,F(L,M)).
\]
These observations yield the following result: 

\begin{thm}\label{textmap}
Let $R$ be an equivariant $S$-algebra, and $L$ and $M$ be left 
$R$-modules.  There is a natural
transformation of graded Mackey functors
\[
\TExt_{R}^{-*}(L,M)\lrarrow
\rfunc{(\pimac R)}{\pimac L}{\pimac M}.
\]
If $L$ is a weak bimodule, this is a map of left
$\pimac R$-modules.  If $M$ is a weak bimodule, it is a map of
right $\pimac R$-modules.
\end{thm}

%

The remainder of this section is devoted to the statements 
and proofs of two results about $R$-modules whose homotopy Mackey functors
are projective or injective as $\pR$-modules.  The behavior of the 
natural maps of Theorems~\ref{ttormap} and~\ref{textmap} for such $R$-modules 
is of particular interest to us.  For this discussion, we denote 
$\pimac R$ by $\pR$.  

\begin{thm}\label{projinjexist}
\indent
\begin{thmlist}
\item If $\mgMP$ is a projective left $\pR$-module, then there
exists a left $R$-module $P$ such that $\pimac P \iso \mgMP$. 
\item If $\mgMQ$ is a projective right $\pR$-module, then there
exists a right $R$-module $Q$ such that $\pimac Q \iso \mgMQ$. 
\item\llabel{pieiii} If $\mgMI$ is an injective left $\pR$-module, then there
exists a left $R$-module $I$ such that $\pimac I \iso \mgMI$.
\end{thmlist}
\end{thm}

\begin{thm}\label{projinjuniq}
Let $L$ and $M$ be left
$R$-modules, and let $N$ be a
right $R$-module.
\begin{thmlist}
\item\llabel{piui} If $\pimac L$ is projective or $\pimac M$ is 
injective as a left $\pR$-module, then the natural map
$\TExt_{R}^{-*}(L,M)\to
\rfunc{\pR}{\pimac L}{\pimac M}$ is an
isomorphism. 
\item If $\pimac M$ is a projective left $\pR$-module or
$\pimac N$ is a projective right $\pR$-module, then the
natural map $\pimac N\rbox[\pR]\pimac M\to
\TTor^{R}_{*}(N,M)$ is an isomorphism. 
\end{thmlist}
\end{thm}

For any two left $R$-modules $L$ and $M$, $\MExt_{R}^{0}(L,M)(G/G)$ is
canonically isomorphic to the abelian group of
maps from $L$ to $M$ in the derived category of left $R$-modules.  
Likewise. for any left $\mgMR$-modules $\mgML$ and $\mgMM$,
$\rfunc[0]{\mgMR}{\mgML}{\mgMM}(G/G)$ is canonically isomorphic to
the abelian group of maps from $\mgML$ to 
$\mgMM$ in the category of left $\mgMR$-modules.  
This implies the following corollary of Theorem \ref{projinjuniq}.

\begin{cor}\label{mapcor}
Let $L$ and $M$ be left $R$-modules.  If $\pimac L$ is projective 
or $\pimac M$ is injective as a left
$\pR$-module, then maps from $L$ to $M$ in the derived
category of left $R$-modules are in one-to-one correspondence with
maps from $\pimac L$ to $\pimac M$ in the category of left
$\pR$-modules. 
\end{cor}

We begin the proof of Theorems~\ref{projinjexist}
and~\ref{projinjuniq} with the following special case.
We state and prove it in the case of left $R$-modules but the
analogous result holds for right $R$-modules.

\begin{lem}\label{spherelem}
Let $X$ be a finite $G$-set and let $\tindex$ be an element of \ROG.
Then there exists a left $R$-module $\cell{\tindex}{X}$ with $\pimac
\cell{\tindex}{X}\iso
\pR\gbxprd\Sigma^{\tindex}\mgMBr^{X}$. Moreover, for
any left $R$-module $M$ and any right $R$-module $N$, the natural maps
\begin{gather*}
(\pimac N)\rbox[\pR](\pimac \cell{\tindex}{X})\lrarrow
\MTor^{R}_{*}(N,\cell{\tindex}{X}),\\
\TExt_{R}^{-*}(\cell{\tindex}{X},M)\lrarrow
\rfunc{\pR}{\pimac \cell{\tindex}{X}}{\pimac M}
\end{gather*}
are isomorphisms.
\end{lem}

\begin{proof}
By taking $\cell{\tindex}{X}=\cell{0}{X}\sma S^{\tindex}$, it suffices
to consider the case when $\tindex=0$.  Then we take
$\cell{}{X}=\cell{0}{X}=\free{R} \Sigma^{\infty}X_{+}$.
Theorem~\ref{pisym} gives us a 
canonical map of $\pR$-bimodules 
\[
\pR \gbxprd \mgMBr^{X}=
\pimac R \gbxprd \pimac[0] (\Sigma^{\infty} X_{+})\lrarrow
\pimac \cell{}{X},
\]
and an easy Spanier--Whitehead duality argument shows that this map is
an isomorphism.  The argument generalizes to show that the
map $\pimac N\gbxprd \mgMBr^{X}\to \pimac(N\sma X_{+})$ is an
isomorphism for all $N$, and the rest of the proof is an easy check of
diagrams.
\end{proof}

We generalize this to other projective modules in the following lemma.

\begin{lem}\label{projlem}
Let $\mgMP$ be a projective left $\pR$-module.  There exists a left
$R$-module $P$ with $\pimac P\iso \mgMP$ such that the natural maps
$\pimac (-) \rbox[\pR]\mgMP\to \TTor^{R}_{*}(-,P)$ and $\TExt_{R}^{-*}(P,-)\to
\rfunc{\pR}{\mgMP}{\pimac (-)}$ are isomorphisms.
\end{lem}

\begin{proof}
Using Proposition~\ref{charproj}, we can find an epimorphism $f\colon
\mgMF\to \mgMP$, where $\mgMF$ is
a direct sum of $\pR$-modules of the form
$\pR\boxprod\Sigma^{\tindex}\mgMBr^{X}$.
Choose a splitting map $g\colon \mgMP\to \mgMF$ for $f$.
By the previous lemma, there is 
a left $R$-module $F$, which is a wedge of modules of the form
$\cell{\tindex}{X}$, such that $\pimac F\iso \mgMF$.  Also, there is  
a self-map $h\colon F\to
F$ in the derived category of left $R$-modules that induces $g\circ f$
on $\pimac $.  Form the left 
$R$-module $P=h^{-1}F$ as the telescope of the self-map $h$.
Then $\pimac P\iso \mgMP$.  For any right $\pR$-module $\mgMN$ and any
right $R$-module $N$, the natural maps
\begin{gather*}
\Colim \mgMN\rbox[\pR]\mgMF\lrarrow \mgMN\rbox[\pR](\Colim \mgMF)\iso 
\mgMN \rbox[\pR] \mgMP,\\
\Colim \TTor^{R}_{*}(N,F)\lrarrow
\TTor^{R}_{*}(N,\Tel F)= 
\TTor^{R}_{*}(N,P)
\end{gather*}
associated to the sequential colimits over the
self-maps $g\circ f$ and $h$ are isomorphisms.  
By Lemma~\ref{spherelem}, the natural
map $\pimac (-)\rbox[\pR]\mgMF\to \TTor^{R}_{*}(-,F)$ is an isomorphism, and
so the natural map $\pimac (-)\rbox[\pR]\mgMP\to \TTor^{R}_{*}(-,P)$ 
is also an isomorphism.  

For any left $R$-module $M$, we have the usual 
short exact sequences of homotopy groups 
\[
0\to \Lim^{1} \pi^{H}_{*+1}\derF_{R}(F,M)
\to \pi^{H}_{*}\derF_{R}(P,M) 
\to \Lim \pi^{H}_{*} \derF_{R}(F,M) \to 0.  
\]
associated to a 
telescope.  
Since $(g\circ f)\circ (g\circ
f)=(g\circ f)$, Lemma~\ref{spherelem} implies that $h\circ h=h$.
It follows that the towers of abelian groups in
question are Mittag-Leffler, and so $\Lim^{1}=0$.  Thus, the natural map 
\[
\TExt_{R}^{-*}(P,M) \to \Lim \TExt_{R}^{-*}(F,M)
\]
is an isomorphism.  For any left $\pR$-module $\mgMM$, the functor 
$\rfunc{\pR}{-}{\mgMM}$ converts colimits to limits, and so the
natural map
$\rfunc{\pR}{\mgMP}{\mgMM}\to \Lim \rfunc{\pR}{\mgMF}{\mgMM}$
is an isomorphism.  By
Lemma~\ref{spherelem} again, 
the natural map $\TExt_{R}^{-*}(F,-)\to
\rfunc{\pR}{\pF}{\pimac (-)}$ is an isomorphism.  The
natural map $\TExt_{R}^{-*}(P,-)\to
\rfunc{\pR}{\mgMP}{\pimac (-)}$ is therefore also an isomorphism.
\end{proof}

If $P'$ is any other left $R$-module with $\pimac P'$
isomorphic to a projective left $\pR$-module $\mgMP$, then $P'$ is
isomorphic in the derived category of left $R$-modules to the left
$R$-module $P$ of the previous lemma.  To see this, note that a 
special case of the
isomorphism $\MExt_{R}^{-*}(P,P')\iso
\rfunc{\pR}{\pimac P}{\pimac P'}$ of the previous lemma 
indicates that there is a one-to-one
correspondence between maps from $P$ to $P'$
in the derived category of left $R$-modules
and 
maps from $\pimac P$ to
$\pimac P'$ in the category of left $\pR$-modules.  Choosing an 
isomorphism $\pimac P\iso \mgMP\iso
\pimac P'$, we obtain a map $P\to P'$ inducing an isomorphism on
homotopy groups.  This proves the following proposition.

\begin{prop}
If $P$ is a left $R$-module such that $\pP=\pimac P$ is a projective
left $\pR$-module, then the natural maps $\pimac (-)\rbox[\pR]\pP\to
\TTor^{R}_{*}(-,P)$ and $\TExt_{R}^{-*}(P,-)\to
\rfunc{\pR}{\pP}{\pimac (-)}$ are isomorphisms.
\end{prop}
 
This gives half of part~(\ref{piui}) of Theorem~\ref{projinjuniq}.
The other half is given by the following lemma.

\begin{lem}\label{ExtHomIso}
If $I$ is a left $R$-module such that $\pI=\pimac I$ is an
injective left $\pR$-module, then the natural map
$\TExt_{R}^{-*}(-,I)\to \rfunc{\pR}{\pimac (-)}{\pI}$ is an isomorphism.
\end{lem}

\begin{proof}
Let $\aC_{I}$ denote the class of left $R$-modules $L$ 
for which the map
$\TExt_{R}^{-*}(L,I)\to \rfunc{\pR}{\pimac L}{\pI}$ is an
isomorphism. Clearly $\aC_{I}$ is closed under arbitrary 
wedge products 
and under suspension by any element $\tindex$ of \ROG.  
Also, an $R$-module $L$ is in $\aC_{I}$ if any $R$-module 
isomorphic to $L$ in the derived category is in $\aC_{I}$.  
By Lemma~\ref{spherelem}, the modules
$\cell{n}{G/H}$ are
in $\aC_{I}$.  Theorem~\ref{homological} indicates that 
the functor $\rfunc{\pR}{-}{\pI}$
is exact.  Thus, if 
\[
\cdots \to \Sigma^{n-1}C 
\to \Sigma^{n}A \to \Sigma^{n}B \to \Sigma^{n}C
\to \Sigma^{n+1}A\to \cdots 
\]
is a cofibration sequence, then applying either of the functors 
$\rfunc{\pR}{\pimac (-)}{\pI}$ or 
$\TExt_{R}^{-*}(-,I)$ to this sequence produces a long exact 
sequence of graded
Mackey functors.  It follows that the cofiber of a map between
modules in $\aC_{I}$ is a module in $\aC_{I}$.  From this, we 
conclude that every left
$R$-module is in $\aC_{I}$.
\end{proof}

The right module parts of Theorems~\ref{projinjexist}
and~\ref{projinjuniq} can be proven by arguments analogous 
to those already given for left $R$-modules.  Thus, to complete 
the proofs of these two results, it suffices to prove 
part~(\ref{pieiii}) of Theorem~\ref{projinjexist}.  The construction 
of the required $R$-modules with injective
homotopy groups is completely analogous to the construction of
Brown--Comenetz dual spectra.  

Let $\mgMI$ be an injective left $\pR$-module.  We define a
contravariant functor $F_{I}$ from the derived category of left
$R$-modules to the category of abelian groups by letting $F_{I}(M)$ 
be the abelian group of maps of left $\pR$-modules from $\pimac M$
to $\mgMI$.  The derived category of left $R$-modules and the functor
$F_{I}$ satisfy the hypotheses for the abstract form of the Brown
representability theorem in Brown~\cite{brown}.  It follows that there
exists a left $R$-module $I$ representing $F_{I}$, i.e., the abelian
group of maps in the derived category of left $R$-modules from $M$ to
$I$ is naturally isomorphic to $F_{I}(M)$.  In particular, letting $M$
range over the left $R$-modules $\cell{\tindex}{X}$, we see that
$\pimac I\iso \mgMI$.  This completes the proof of 
Theorem~\ref{projinjexist}.  Note that Lemma \ref{ExtHomIso} 
implies that the $R$-module $I$ is unique up to isomorphism in 
the derived category of left $R$-modules.  


\section{The Construction of the Spectral Sequences}
\label{secss}

The two spectral sequences described in the introduction 
are constructed in this section.  Throughout this construction, 
$R$ is a fixed equivariant $S$-algebra and $\pR=\pimac R$.  The
results in the previous section allow us to construct 
``resolutions'' of an $R$-module $M$ 
in the derived category of $R$-modules corresponding to 
projective and injective $\pR$-module resolutions of $\pimac M$.  
These resolutions are the equivariant generalization of the 
resolutions constructed in Section IV.5 of \cite{EKMM}.  

The following definition formalizes the relationship between our 
topological resolutions of $M$ and 
algebraic resolutions of $\pimac M$.  

\begin{defn}\label{defres}
Let $M$ be a left $R$-module, and let $\pM=\pimac M$.
A projective topological resolution of $M$ consists of  
collections of left $R$-modules $M_{s}$ and $P_{s}$ together 
with cofiber sequences
\[
\Sigma^{s}P_{s}\overto{j_{s}}
M_{s}\overto{i_{s+1}} 
M_{s+1}\overto{k_{s}} 
\Sigma^{s+1} P_{s},
\]
in the derived category of left $R$-modules
for all $s\geq 0$.  These objects must satisfy the conditions that 
$M_{0}=M$, each $\pimac P_{s}$ is a
projective left $\pR$-module, and $\pimac  j_{s}$ is an epimorphism.

A projective topological resolution $(M_{s},P_{s})$ of $M$ is 
compatible with a projective resolution $\mac{P}_{*,*}$ of $\pM$ if there are 
isomorphisms $\mac{P}_{s,*}\to \pimac P_{s}$ under which 
$\pimac j_{0}$ coincides with the augmentation $\mac{P}_{0,*}\to
\pM$ and $\pimac (k_{s}\circ j_{s+1})$ coincides with the 
suspension $\Sigma^{s+1} d_{s+1}$ of the differential $d_{s+1}\colon \mac{P}_{s+1,*}\to \mac{P}_{s,*}$.  

An injective topological resolution of $M$ consists of collections of
left $R$-modules $M^{s}$ and $I^{s}$ together with fiber sequences 
\[
\Omega^{s+1}I^{s}\overto{k_{s}}
M^{s+1}\overto{i^{s+1}} 
M^{s}\overto{j^{s}} 
\Omega^{s} I^{s},
\]
in the derived
category of left $R$-modules
for all $s\geq 0$.  These must satisfy the conditions that $M^{0}=M$, 
each $\pimac I^{s}$ is
an injective left $\pR$-module, and $\pimac j^{s}$ is a
monomorphism.  

An injective topological resolution $(M^{s},I^{s})$ of $M$ is 
compatible with a given injective resolution $\mac{I}_{*}^{*}$ 
of $\pM$ if there are isomorphisms $\mac{I}_{*}^{s}\to
\pimac I^{s}$ under which $\pimac j^{0}$ coincides with the  
augmentation $\pM\to \mac{I}_{*}^{0}$ and $\pimac (j^{s+1}\circ k^{s})$ 
coincides with the desuspension $\Omega^{s+1} d^{s}$ of the 
differential $d^{s}\colon \mac{I}_{*}^{s}\to \mac{I}_{*}^{s+1}$.
\end{defn}

Projective topological resolutions of right $R$-modules are defined 
analogously.  In the projective context, the modules $M_{s}$ are
analogous to the quotients $X/X^{s-1}$ for a nice filtration
$X_{0}\subset X_{1}\subset X_{2}\subset \cdots \subset X$ of a space
$X$. As discussed in Boardman~\cite[12.5ff]{CCvrgtSS}, the cohomological
spectral sequence constructed from the cofiber sequences
\[
X/X^{s-1}\lrarrow X/X^{s}\lrarrow X^{s}/X^{s-1}
\]
is isomorphic to the more usual spectral sequence constructed from
the cofiber sequences
\[
X^{s-1}\lrarrow X^{s}\lrarrow X^{s}/X^{s-1},
\]
but has better structural properties.  In the derived category
context, Verdier's octahedral axiom converts ``filtrations'' (like the
sequence $X^{0}\to X^{1}\to X^{2}\to\cdots $ or the sequence $\bar
M_{0}\to \bar M_{1}\to \bar M_{2}\to\cdots $ used in the next
section) into ``resolutions'' (like $X\to X/X^{0}\to X/X^{1}\to
\cdots$ or $M_{0}\to M_{1}\to M_{2}\to \cdots$) and vice-versa, but
not canonically.

As noted below, the requirement in the definition of projective topological resolution
that each $\pimac j_{s}$ is an epimorphism ensures that projective 
topological resolutions induce projective algebraic resolutions.  
An analogous observation applies to injective
topological resolutions.

\begin{prop}
If $(M_{s},P_{s})$ is a projective topological resolution of $M$, then
$\pimac P_{s}$ is a complex of $\pR$-modules with differential 
$d_{s+1}=\Sigma^{-(s+1)}\pimac(k_{s}\circ j_{s+1})$ and is a projective 
$\pR$-module resolution 
of $\pimac M$ with augmentation $\pimac j_{0}$.
\end{prop}

It is somewhat less obvious but quite important for us that there is a 
topological resolution corresponding to every projective or injective 
algebraic resolution.  

\begin{lem}
Let $M$ be a left $R$-module.  Then every projective $\pR$-module 
resolution of
$\pM=\pimac M$ has a compatible projective topological resolution.
Also, every injective $\pR$-module resolution of $\pM$ has 
a compatible injective topological resolution.
\end{lem}

\begin{proof}
We treat the projective case in detail; the injective case is entirely
analogous.  Let $\mac{P}_{*,*}$ be a projective resolution of $\pM$.
First apply Theorem~\ref{projinjexist} to choose left $R$-modules
$P_{s}$ with $\pimac P_{s}\iso \mac{P}_{s,*}$.  By
Corollary~\ref{mapcor}, there is a unique map 
$j_{0}\colon P_{0}\to M_{0} = M$ in the derived category that corresponds 
on passage to homotopy Mackey functors to the
augmentation $\mac{P}_{0,*}\to \pM$.  Take $M_{1}$ to be the cofiber
of $j_{0}$, and let $i_{1}$ and $k_{0}$ be the appropriate maps from the 
resulting cofiber sequence.  Now assume by
induction that the resolution has been constructed up to $M_{s}$ and
that $\pimac k_{s-1}$ is injective and induces an isomorphism of
$\pimac M_{s}$ with the submodule $\Sigma^{s}\Image(d_{s})$ of
$\Sigma^{s}\mac{P}_{s-1,*}$.  
Then, by Corollary~\ref{mapcor}, there exists a map
$j_{s}\colon \Sigma^{s}P_{s}\to M_{s}$ whose induced map on homotopy
Mackey functors corresponds to the map 
$d_{s}\colon \mac{P}_{s,*}\to \Image(d_{s})$ under
the chosen isomorphisms.  
Take $M_{s+1}$ to be the cofiber of
$j_{s}$, and let $i_{s+1}$ and $k_{s}$ be the induced maps.  Since the
map $j_{s}$ induces an epimorphism on homotopy Mackey functors, the long
exact sequence associated to this cofiber sequence is short exact.  It follows 
that the map on homotopy Mackey functors induced by $k_{s}$ provides 
an isomorphism between 
$\pimac M_{s+1}$ and $\Sigma^{s+1}\Image(d_{s+1})$.  This completes the
induction.
\end{proof}

Now we are ready to construct the spectral sequences.  Let $M$ be a left
$R$-module and $N$ be a right $R$-module.  Let $\mac{P}_{*,*}$ be a
projective resolution of $\pM=\pimac M$, and let
$(M_{s},P_{s})$ be a projective topological resolution of $M$ compatible with
$\mac{P}_{*,*}$.  Extend the
collection of cofiber sequences to negative $s$ by setting $P_{s}=*$ 
and $M_{s}=M$ for $s<0$ (with $i_{s+1}=\id$ and $j_{s}$ and
$k_{s}$ the trivial map).  Applying the functor $\TTor_{*}^{R}(N,-)$ to
these cofiber sequences yields a homologically graded exact couple with 
\begin{equation}\label{lefttor}
\begin{split}
\mac{D}_{s,\tindex} &= \TTor^{R}_{s+\tindex}(N,M_{s})\\
\mac{E}_{s,\tindex} &= \TTor^{R}_{\tindex}(N,P_{s}).
\end{split}
\end{equation}
The maps 
\[
\spreaddiagramcolumns{-2pc}
{\diagram
\relax\cdots\rrto^{i}&&\relax\mac{D}_{s-1,*}\rrto^{i}
&&\relax\mac{D}_{s,*}\rrto^{i}\dlto^{k}
&&\relax\mac{D}_{s+1,*}\rrto^{i}\dlto^{k}
&&\relax\mac{D}_{s+2,*}\rrto^{i}\dlto^{k}
&&\relax\cdots\\
&\relax\cdots&&\relax\mac{E}_{s-1,*}\ulto^{j}
&&\relax\mac{E}_{s,*}\ulto^{j}
&&\relax\mac{E}_{s+1,*}\ulto^{j}
&&\relax\cdots
\enddiagram}
\]
in this exact couple come from the analogously named maps in 
the topological resolution.  Both $i$ and $j$ preserve the 
total degree $s+\tindex$, and $k$
lowers it by one.  Alternatively, given
a projective resolution $\mac{Q}_{*,*}$ of $\pN=\pimac N$ and a
compatible projective topological resolution $(N_{s},Q_{s})$ of $N$, setting
\begin{equation}\label{righttor}
\begin{split}
\mac{D}_{s,\tindex} &= \TTor^{R}_{s+\tindex}(N_{s},M)\\
\mac{E}_{s,\tindex} &= \TTor^{R}_{\tindex}(Q_{s},M)
\end{split}
\end{equation}
gives another homologically graded exact couple of exactly the same form. 

For the $\Ext$ spectral sequence, let $L$ and $M$ be left $R$-modules.  
Also, let $\mac{O}_{*,*}$ be a
projective resolution of $\pL=\pimac L$ and let $(L_{s},O_{s})$ be
a compatible projective topological resolution.  Again set $L_{s}=L$,
$O_{s}=*$ for $s<0$.  Apply the functor $\TExt_{R}^{*}(-,M)$ to the
cofiber sequences relating $O_{s}$ and $L_{s}$, and let
\begin{equation}\label{leftext}
\begin{split}
\mac{D}^{s,\tindex} &= \TExt_{R}^{s+\tindex}(L_{s},M)\\
\mac{E}^{s,\tindex} &= \TExt_{R}^{\tindex}(O_{s},M).
\end{split}
\end{equation}
This yields a cohomologically graded exact couple of the form
\[
\spreaddiagramcolumns{-2pc}
{\diagram
\relax\cdots\rrto^{i}&&\relax\mac{D}^{s+2,*}\rrto^{i}
&&\relax\mac{D}^{s+1,*}\rrto^{i}\dlto^{j}
&&\relax\mac{D}^{s,*}\rrto^{i}\dlto^{j}
&&\relax\mac{D}^{s-1,*}\rrto^{i}\dlto^{j}
&&\relax\cdots\\
&\relax\cdots&&\relax\mac{E}^{s+1,*}\ulto^{k}
&&\relax\mac{E}^{s,*}\ulto^{k}
&&\relax\mac{E}^{s-1,*}\ulto^{k}
&&\relax\cdots\ .
\enddiagram}
\]
In this exact couple, the maps $i$ and $j$ preserve the 
total cohomological degree $s+\tindex$,
and $k$ raises it by one.  Alternatively, let $\mac{I}_{*}^{*}$ be 
an injective resolution of $\pM$ and 
$(M^{s},I^{s})$ be a compatible injective topological
resolution of $M$.  Extend the fiber sequences of the topological 
resolution by setting $I^{s}=*$ and
$M^{s}=M$ for $s<0$.  Applying $\TExt_{R}^{*}(L,-)$ to the fiber
sequences relating $M^{s}$ and $I^{s}$ and setting
\begin{equation}\label{rightext}
\begin{split}
\mac{D}^{s,\tindex} &= \TExt_{R}^{s+\tindex}(L,M^{s})\\
\mac{E}^{s,\tindex} &= \TExt_{R}^{\tindex}(L,I^{s})
\end{split}
\end{equation}
gives a cohomologically graded exact couple of the same form as above.


These exact couples lead to spectral sequences in the usual way.
These spectral sequences are clearly natural in the
unresolved variable.  Theorem~\ref{projinjuniq} identifies the
$E^{1}$- and $E^{2}$-terms.  When the unresolved variable is a weak
bimodule, Theorems~\ref{ttormap} and~\ref{textmap} imply that the exact
couple is an exact couple of $\pR$ or $\pcR$-modules.  The resulting spectral
sequence is therefore a spectral sequence of $\pR$ or $\pcR$-modules.  These 
assertions are summarized in the  following result:

\begin{thm}
Let $L$ and $M$ be left $R$-modules and $N$ be a right $R$-module.
\begin{thmlist}
\item The spectral sequence derived from the exact
couple~(\ref{lefttor}) has $E^{1}$ complex canonically isomorphic to
$\pN\rbox[\pR]\mac{P}_{*,*}$ and $E^{2}_{s,\tindex}$-term canonically
isomorphic to $\MTor^{\pR}_{s,\tindex}(\pN,\pM)$.  The spectral
sequence is natural in $N$.  If $N$ is a weak bimodule, this is a
spectral sequence of left $\pR$-modules.
\item The spectral sequence derived from the exact
couple~(\ref{righttor}) has
$E^{1}$ complex canonically isomorphic to
$\mac{Q}_{*,*}\rbox[\pR]\pM$ and $E^{2}_{s,\tindex}$-term 
canonically isomorphic to 
$\MTor^{\pR}_{s,\tindex}(\pN,\pM)$.  The spectral sequence is
natural in $M$.  If $M$ is a weak bimodule, this is a spectral
sequence of right $\pR$-modules.
\item The spectral sequence derived from the exact
couple~(\ref{leftext}) has 
$E_{1}$ complex canonically isomorphic to
$\rfunc{\pR}{\mac{O}_{*,*}}{\pM}$ and
$E_{2}^{s,\tindex}$-term canonically
isomorphic to 
$\MExt_{\pR}^{s,\tindex}(\pL,\pM)$.  The spectral sequence is
natural in $M$.  If $M$ is a weak bimodule, this is a spectral
sequence of right $\pcR$-modules.
\item The spectral sequence derived from the exact
couple~(\ref{rightext}) has $E_{1}$ complex canonically isomorphic to
$\rfunc{\pR}{\pL}{\mac{I}^{*}_{*}}$ and
$E_{2}^{s,\tindex}$-term canonically
isomorphic to 
$\MExt_{\pR}^{s,\tindex}(\pL,\pM)$.  The spectral sequence is
natural in $L$.  If $L$ is a weak bimodule, this is a spectral
sequence of left $\pcR$-modules.
\end{thmlist}
\end{thm}

This theorem leaves unresolved the question of whether our spectral 
sequences are natural in the resolved variable.  It is also not obvious 
that these spectral sequences are independent of the resolution 
chosen to form them.  To settle these questions, we prove the
following result in the next section.

\begin{thm}\label{siso}
The identity map on $\MTor^{\pR}_{*,*}(\pN,\pM)$ induces an 
isomorphism between the spectral
sequences derived from~(\ref{lefttor}) and~(\ref{righttor}).  Similarly, 
the identity map on $\MExt_{\pR}^{*,*}(\pL,\pM)$ induces an isomorphism 
between the spectral
sequences derived from~(\ref{leftext}) and~(\ref{rightext}).
\end{thm}

The issue of convergence for these spectral sequences must 
still be discussed.  Since limits and colimits of Mackey
functors are formed object-wise, convergence of spectral sequences of
Mackey functors works just like convergence of spectral sequences of
abelian groups.  The spectral sequences derived
from~(\ref{lefttor}) and~(\ref{righttor}) are homological right
half-plane spectral sequences.  Thus, by Boardman \cite[6.1]{CCvrgtSS}, 
to prove 
that they converge strongly to
\[
\TTor^{R}_{\tindex}(N,M) = \Lim_{s} \mac{D}_{-s,s+\tindex},
\]
it suffices to show that $\Colim_{s} \mac{D}_{s,-s+\tindex}=0$. In the 
case of~(\ref{lefttor}), this amounts to showing that $\Colim_{s}
\TTor^{R}_{\tindex}(N,M_{s})=0$.  For this, consider the cofiber
sequence 
\[
\bigvee M_{s} \lrarrow \bigvee M_{s}\lrarrow \Tel M_{s}
\]
defining the telescope of a sequence of maps in the derived category.
The canonical map $\Colim \pimac M_{s}\to \pimac (\Tel M_{s})$
is an isomorphism, and it follows that $\Tel M_{s}$ is trivial.  Since
the derived smash product over $R$ preserves cofiber sequences, 
$\Tel \dersma^{R}({N},{M_{s}})$ is isomorphic to
$\dersma^{R}({N},{\Tel M_{s}})$ and is therefore trivial.  In
particular, $\Colim \TTor^{R}_{*}(N,M_{s})=0$.  Since the edge
homomorphism in this case is induced by $j_{0}$, we obtain the
following result. 

\begin{thm}
The spectral sequence derived from the exact couple~(\ref{lefttor})
converges strongly to $\TTor^{R}_{*}(N,M)$.  Its edge homomorphism is the
canonical map $\pN\rbox[\pR]\pM\to \TTor^{R}_{*}(N,M)$.
\end{thm}

For cohomologically graded right half-plane spectral sequences like those derived
from~(\ref{leftext}) and~(\ref{rightext}), conditional convergence to 
\[
\TExt_{R}^{\tindex}(L,M) = \Colim_{s} \mac{D}^{-s,s+\tindex}
\]
is defined (see, for example, Boardman~\cite[5.10]{CCvrgtSS}) to mean that 
\[
\Lim_{s} \mac{D}^{s,-s+\tindex}=0\quad\text{and}\quad
\Lim^{1}_{s} \mac{D}^{s,-s+\tindex}=0.
\]
In the context of~(\ref{rightext}),
$\mac{D}^{s,-s+\tindex}=\TExt_{R}^{\tindex}(L,M^{s})$.  The limit and 
$\Lim^{1}$ term that must vanish are defined by an exact sequence
\[
0\to
\Lim_{s} \mac{D}^{s,-s+\tindex}\to
\prod \TExt_{R}^{\tindex}(L,M^{s})\to
\prod \TExt_{R}^{\tindex}(L,M^{s})\to
\Lim^{1}_{s} \mac{D}^{s,-s+\tindex}\to 0.
\]
Consider the fiber sequence 
\[
\Mic M^{s}\lrarrow \prod M^{s}\lrarrow \prod M^{s}
\]
defining the microscope of a sequence of maps in the derived
category.  Each map $i^{s+1}\colon M^{s+1}\to M^{s}$ induces the zero
map on homotopy Mackey functors, and so $\Mic M^{s}$ is trivial.  The derived
function spectrum functor $\derF_{R}({L},{-})$ preserves fiber
sequences, and so
\[
\derF_{R}({L},{\Mic M^{s}})
\lrarrow \prod \derF_{R}({L},{M^{s}})
\lrarrow \prod \derF_{R}({L},{M^{s}})
\]
is a fiber sequence.  Since $\Mic M^{s}$ is trivial, so is
$\derF_{R}({L},{\Mic M^{s}})$.  The induced long exact sequence on
homotopy Mackey functors now gives conditional convergence.  Since 
the edge homomorphism in this context is induced by $j^{0}$, we obtain 
the following result. 

\begin{thm}
The spectral sequence derived from the exact couple
of~(\ref{rightext}) converges conditionally to $\TExt_{R}^{*}(L,M)$.  Its 
edge homomorphism is the canonical map $\TExt_{R}^{*}(L,M)\to
\rfunc{\pR}{\pL}{\pM}$.
\end{thm}

This completes the construction of the Hyper-Tor and Hyper-Ext spectral
sequences described in the introduction.  


\section{Uniqueness, Naturality, and The Yoneda Pairing}
\label{secyon}

In the previous section, we constructed a pair of Hyper-Tor spectral
sequences and a pair of Hyper-Ext spectral sequences.  This section 
contains the proof of Theorem~\ref{siso}, which asserts that the pair 
of Hyper-Tor spectral sequences are isomorphic and so are the pair of 
Hyper-Ext spectral sequences.  The technical work needed to prove this 
result also suffices to construct construct the Yoneda pairing of 
Hyper-Ext spectral sequences mentioned in the introduction.  This 
construction is described in Theorem~\ref{ssypair}. 

The formation of our spectral sequences in the last section did not 
require any constructions more complicated than cofiber and fiber 
sequences, telescopes, and microscopes.  However, here we need more 
complicated homotopy colimits and limits.  Because of this, the 
constructions described in this section must be carried out in a point 
set category rather than the associated derived category.  Nevertheless, 
these constructions are of such a general nature that they can be 
carried out in any of the modern categories of equivariant spectra 
(i.e., \cite{EKMM,GOrtho,gss}).  

For the remainder of the section, fix left $R$-modules
$L$ and $M$ and a right $R$-module $N$.  Without loss of generality, 
it can be assumed that each of these objects is cofibrant and fibrant 
in the appropriate module category.  Also, fix projective
topological resolutions $(L_{s},O_{s})$, $(M_{s},P_{s})$, and
$(N_{s},Q_{s})$ of $L$, $M$, and $N$, respectively, and an injective
topological resolution $(M^{s},I^{s})$ of $M$.  We prove 
Theorem~\ref{siso} by constructing isomorphisms between the pairs 
of spectral sequences derived from these specific resolutions.  

As noted in the remarks preceding the definition of a projective 
topological resolution, the sequence of $R$-modules 
\[
M = M_{0} \lrarrow M_{1} \lrarrow \cdots \lrarrow M_{s} \lrarrow \cdots
\]
is analogous to the sequence of quotients of $M$ 
by a sequence of progressively larger submodules of $M$,
\[
M = M/\bar M_{-1} \lrarrow M/\bar M_{0} \lrarrow \cdots
\lrarrow M/\bar M_{s-1} \lrarrow \cdots .
\]
To form the 
homotopy limits and colimits needed for the proof of Theorem~\ref{siso}, 
we must ``reconstruct'' these missing submodules.  Specifically, we
construct a sequence of (point-set level) cofibrations of $R$-modules
\[
*= \bar M_{-1}\lrarrow \bar M_{0} \overto{\cofm{1}} \bar M_{1} \lrarrow  
\cdots \lrarrow \bar M_{s-1} \overto{\cofm{s}} \bar M_{s} \lrarrow  \cdots
\]
together with compatible point-set level $R$-module maps
$\incm{s-1}\colon \bar M_{s-1}\to M$ whose behavior in the derived
category is what one would expect from an appropriate filtration of
$M$ with quotients $M_{s}$.  In particular, for each $s$, we choose
isomorphisms $C(\incm{s-1}) \iso M_{s}$ and $C(\cofm{s}) \iso
\Sigma^{s} P_{s}$ in the derived category that are
compatible in the following sense: They provide an isomorphism in the
derived category between the induced maps of homotopy cofibers
\begin{equation}\label{cofiber}
C(\cofm{s}) \lrarrow C(\incm{s-1}) \lrarrow C(\incm{s}) \lrarrow 
C(\Sigma \cofm{s})\simeq \Sigma C(\cofm{s})
\end{equation}
and the cofiber sequence 
\begin{equation}\label{cofiberorig}
\Sigma^{s} P_{s}\lrarrow M_{s}\lrarrow M_{s+1}\lrarrow \Sigma^{s+1}P_{s}
\end{equation}
that is a part of our chosen projective topological resolution.
Here and for the rest of the section, we understand the point-set
model for homotopy
cofibers to be formed using the usual cone construction; that is, 
\[
C(\cofm{s}) = \bar M_{s}\cup_{\bar M_{s-1}}(\bar M_{s-1}\sma I_{+})
\cup_{\bar M_{s-1}}*
\]
and
\[
C(\incm{s-1}) = M\cup_{\bar M_{s-1}}(\bar M_{s-1}\sma I_{+})
\cup_{\bar M_{s-1}}*.
\]
Note that since the model categories in \cite{EKMM,GOrtho,gss} are
simplicial and the objects $\bar M_{s-1}$ and $M$ are cofibrant, 
these are cofibrant objects and 
are cofibers in the sense of Quillen.  Moreover, sequence (\ref{cofiber})
represents a cofiber sequence in the derived category.  The apparent
shift in indexing in the cofiber sequence
\[
\bar M_{s-1}\overto{\incm{s-1}}M\lrarrow M_{s}\lrarrow \Sigma \bar M_{s-1}
\]
($M_{s}$ is the cofiber of $\bar M_{s-1}$)
is for consistency with the grading conventions of
Boardman~\cite[\S12]{CCvrgtSS}.  This shift also leads to cleaner
formulas in the work below. 

The ``reconstruction'' of the modules $\bar M_{s}$ essentially amounts to a
point-set refinement of Verdier's octahedral axiom.  The standard
proof of this axiom is actually strong enough to provide the refinement
we need.  Let $\bar M_{-1}=*$.  Next choose a cofibrant $R$-module
$\bar M_{0}$ weakly equivalent to $P_{0}$ and a
point-set map $\incm{0}\colon \bar M_{0}\to M$ of $R$-modules such that
the induced cofiber sequence 
\[
\bar M_{0}\lrarrow M \lrarrow C(\incm{0}) \lrarrow \Sigma \bar M_{0}
\]
is isomorphic in the derived category to the given cofiber sequence 
\[
P_{0}\lrarrow M \lrarrow M_{1}\lrarrow \Sigma P_{0}
\]
by an isomorphism that is the identity on $M$. 
For the inductive step, assume that the sequence of 
cofibrations and compatible maps into $M$ has been constructed up 
to the $R$-module $\bar M_{s}$.  Choose a fibrant
$R$-module $X$ and a point-set level $R$-module map
$C(\incm{s}) \to X$ such that the induced cofiber sequence is isomorphic in 
the derived category to the given cofiber sequence 
\[
M_{s+1}\lrarrow M_{s+2}\lrarrow \Sigma^{s+2} P_{s+1} \to \Sigma M_{s+1}
\]
by an isomorphism that restricts to the previously chosen isomorphism
$C(\incm{s})\to M_{s+1}$.  Let $F$ be the homotopy fiber of 
the composite 
$M\to C(\incm{s}) \to X$
defined using the usual path space construction,
\[
F=M\times_{X}X^{I_{+}}\times_{X}*.
\]
Since $M$ and $X$ are fibrant, this represents the fiber in the sense
of Quillen.  By construction, the composite 
$\bar M_{s} \to M \to X$
factors through the composite $\bar M_{s}\to C\bar M_{s}\to C(\incm{s})$.  
This factorization provides a null homotopy of the map from 
$\bar M_{s}$ to $X$ and so a lift $\bar M_{s}\to F$ of
the map $\bar M_{s}\to M$.  Choose $\bar M_{s+1}$ by factoring the
map $\bar M_{s}\to F$ as a cofibration $\cofm{s+1}\colon \bar M_{s}\to
\bar M_{s+1}$ followed by an acyclic fibration $\bar
M_{s+1}\overto{\sim}F$.  Let $\incm{s+1}\colon \bar M_{s+1}\to M$ be
the composite map $\bar M_{s+1}\to F\to M$.  The objects and maps we 
have chosen fit into the diagram 
\[
\xymatrix{
\relax\bar M_{s}\ar[r]^{\incm{s}}\ar[d]_{\cofm{s+1}}
&M\ar[r]\ar @{=} [d]&C(\incm{s})\ar[d]\\
\relax\bar M_{s+1}\ar[r]^{\incm{s+1}}\ar @{-->} [d]_{\sim}
&M\ar[r]\ar @{-->} [dr]&C(\incm{s+1})\ar @{-->} [d]^{\sim}\\
F\ar @{-->} [ur]&&X,
}
\]
which commutes on the point-set level.  
The map $C(\cofm{s+1})\to X$ is a weak equivalence by \cite[I.6.4]{EKMM},
\cite[3.5.(vi)]{GOrtho}, or \cite[5.8]{gss}.  Thus, this procedure constructs
an isomorphism in the derived category between cofiber
sequence~(\ref{cofiber}) and cofiber sequence~(\ref{cofiberorig})
that restricts to the previously chosen isomorphism $C(\incm{s})\to
M_{s+1}$.  This completes our construction, which is described formally 
in the following proposition.

\begin{prop}
Let $(M_{s},P_{s})$ be a projective topological resolution of an 
$R$-module $M$.  Then there is a sequence 
\[
*= \bar M_{-1}\lrarrow \bar M_{0} \overto{\cofm{1}} \bar M_{1} \lrarrow  
\cdots \lrarrow \bar M_{s-1} \overto{\cofm{s}} \bar M_{s} \lrarrow  \cdots
\]
of cofibrations of $R$-modules together with $R$-module maps
$\incm{s-1}\colon \bar M_{s-1}\to M$, compatible on the point-set level, 
such that $C(\incm{s-1}) \iso M_{s}$ and $C(\cofm{s}) \iso
\Sigma^{s} P_{s}$ in the derived category.  Moreover, these isomorphisms 
are compatible in the sense that they provide an isomorphism in the 
derived category between the cofiber sequences 
\[
C(\cofm{s}) \lrarrow C(\incm{s-1}) \lrarrow C(\incm{s}) \lrarrow 
C(\Sigma \cofm{s})\simeq \Sigma C(\cofm{s})
\]
and 
\[
\Sigma^{s} P_{s}\lrarrow M_{s}\lrarrow M_{s+1}\lrarrow \Sigma^{s+1}P_{s}.
\]
\end{prop}

In this same manner, choose sequences of cofibrations 
\[
*=L_{-1}\lrarrow \bar L_{0}\overto{\cofl{1}} \bar L_{1}\lrarrow \cdots
\lrarrow \bar L_{s-1}\overto{\cofl{s}} \bar L_{s}\lrarrow \cdots 
\]
and 
\[
*=N_{-1}\lrarrow \bar N_{0}\overto{\cofn{1}} \bar N_{1}\lrarrow \cdots
\lrarrow \bar N_{s-1}\overto{\cofn{s}} \bar N_{s}\lrarrow \cdots 
\]
together with compatible maps $\incl{s}\colon \bar L_{s} \to L$ and
$\incn{s}\colon \bar N_{s}\to N$ 
consistent with the chosen projective resolutions of $L$ and $N$.  
When it is desirable to indicate which of the objects $L$, $M$, and $N$ 
is associated to a particular map $\cofm{s}$ or $\incm{s}$, the symbols 
$\cofM{s}$, $\incM{s}$, etc., are employed.     

An analogous sequence of fibrations 
\[
\cdots \lrarrow M^{s} \overto{\fibm{s-1}} M^{s-1}
\lrarrow \cdots \lrarrow \bar M^{0}
\overto{\fibm{1}} \bar M^{0}\lrarrow M^{-1}=*
\]
together with compatible maps $\prom{s}\colon M \to M^{s}$ can be
constructed using the ``Eck\-mann--Hilton'' dual of the argument above
(that is, reverse the direction of all of the arrows, and switch 
cofibrations and fibrations, colimits and limits, and $(-)\sma I_{+}$ and
$(-)^{I_{+}}$).

To make use of the $R$-modules just constructed, we must 
introduce a family of very simple categories.

\begin{defn}
Let $\dcat{}$ denote the category which has as objects the ordered
pairs of natural numbers $(s,t)$ and as maps a
unique map $(s,t)\to (s',t')$ whenever $s\leq s'$ and $t\leq t'$.
Let $\dcat{n}$ denote the full subcategory of $\dcat{}$ consisting
of objects $(s,t)$ with $s+t\leq n$
\end{defn}

The (point-set) smash products $\bar N_{s} \sma_{R} \bar
M_{t}$ may be regarded as a functor from $\dcat{}$ to equivariant $S$-modules,
orthogonal spectra, or symmetric spectra, as appropriate.  Likewise,
the function spectra $F_{R}(\bar L_{s},\bar M^{t})$ form a contravariant
functor from $\dcat{}$ to equivariant $S$-modules, orthogonal spectra,
or symmetric spectra.  Our main tools in this section are homotopy
limits and colimits of these functors.  Our argument requires constructions 
of homotopy limits and colimits that are functorial on the point-set level.  
Any functorial constructions should be adequate.  However, at some points 
in the argument, we assume that the ``bar construction'' models of homotopy limits and colimits (described, 
for example, in \cite[X\S3]{EKMM}) are used in order to fill in certain details. 

\begin{defn}\label{HoLimDef}
Let $T_{n}=\Hocolim_{\dcat{n}}(\bar N_{s}\sma_{R}\bar M_{t})$ and
$U_{n}=\Holim_{\dcat{n}}F_{R}(\bar L_{s},\bar M^{t})$.
There are canonical maps  
\[
g_{n+1}\colon T_{n}\lrarrow T_{n+1}
\qquad\text{and}\qquad
g^{n+1}\colon U^{n+1}\lrarrow U^{n}.   
\]
induced by the inclusion of categories $\dcat{n}\to
\dcat{n+1}$.  Similarly, there are maps 
\[
h_{n}\colon T_{n} \to N\sma_{R} M 
\qquad\text{and}\qquad 
h^{n}\colon F_{R}(L,M) \to U^{n}
\]
induced by the maps $\bar M_{s}\to M$ and the
analogous maps for the $\bar N_{s}$, $\bar L_{s}$, and $\bar M^{s}$
sequences.
\end{defn}

There are also maps 
\begin{gather*}
T_{n}\lrarrow N \sma_{R}\bar M_{n} \qquad  T_{n} \lrarrow \bar N_{n}\sma_{R}M, \\
F_{R}(\bar L_{n},M)\lrarrow U^{n} \qquad  
F_{R}(L,\bar M^{n})\lrarrow U^{n}
\end{gather*} 
induced by the maps $\bar M_{s}\to \bar M_{s+1}$
and the analogous maps for the $\bar N_{s}$, $\bar L_{s}$, and $\bar
M^{s}$ sequences.  Because we are using a functorial construction of
homotopy colimits and limits, these 
are point-set level maps.  Moreover, the diagrams  
\[
\diagram
T_{n}\rto\dto_{g_{n+1}}
&N \sma_{R}\bar M_{n}\dto^{\id \sma \,\cofM{n+1}}
&&T_{n}\rto\dto_{g_{n+1}}
&\bar N_{n}\sma_{R}M\dto\dto^{\cofN{n+1} \sma \,\id}\\
T_{n+1}\rto&N \sma_{R}\bar M_{n+1}
&&T_{n+1}\rto&\bar N_{n+1}\sma_{R}M,
\enddiagram
\]
and the analogous diagrams for the homotopy limits $U^{n}$ commute on
the point-set level.

Since the diagrams above commute on the point-set level, they induce 
canonical maps 
\[
{\diagram
C(g_{s+1})\rto\dto
&C(h_{s})\rto\dto
&C(h_{s+1})\rto\dto
&\Sigma C(g_{s+1})\dto\\
C(\id_{N}\sma \cofM{s+1})\rto
&C(\id_{N}\sma \incM{s})\rto
&C(\id_{N}\sma \incM{s+1})\rto
&\Sigma C(\id_{N}\sma \cofM{s+1})
\enddiagram}
\]
and
\[
{\diagram
C(g_{s+1})\rto\dto
&C(h_{s})\rto\dto
&C(h_{s+1})\rto\dto
&\Sigma C(g_{s+1})\dto\\
C(\cofN{s+1}\sma \id_{M})\rto
&C(\incN{s}\sma \id_{M})\rto
&C(\incN{s+1}\sma \id_{M})\rto
&\Sigma C(\cofN{s+1}\sma \id_{M})
\enddiagram}
\]
of cofiber sequences in the derived category.  The constructions 
$N\sma_{R}(-)$ and $(-)\sma_{R}M$ preserve cofiber sequences.  Thus, 
the isomorphisms characterizing the terms in cofiber
sequence~(\ref{cofiber}) can be used to identify the bottom cofiber
sequences in the two diagrams  
above with those defining the exact couples~(\ref{lefttor})
and~(\ref{righttor}), respectively. 
Setting
\begin{equation}\label{bitor}
\begin{split}
\mac{D}_{s,\tindex} &= \pimac[s+\tindex]C(h_{s})\\
\mac{E}_{s,\tindex} &= \pimac[s+\tindex]C(g_{s})
\end{split}
\end{equation}
gives a third homologically graded exact couple.  Moreover, the 
commuting diagrams above provide maps of exact 
couples from (\ref{bitor}) to both (\ref{lefttor})
and~(\ref{righttor}).  

The following result about these three exact couples implies the part
of Theorem~\ref{siso} applicable to exact couples (\ref{lefttor})
and~(\ref{righttor}). 

\begin{thm}\label{thmbitor}
The spectral sequence derived from~(\ref{bitor}) has its $E^{1}$
complex canonically isomorphic to the total complex of
$\mac{Q}_{*,*}\rbox[\pR]\mac{P}_{*,*}$.  Moreover, under the canonical
isomorphisms, the maps of this $E^{1}$
complex to the $E^{1}$ complexes of (\ref{lefttor}) and
(\ref{righttor}) coincide with the augmentations
\[
\mac{Q}_{*,*}\rbox[\pR]\mac{P}_{*,*}\lrarrow \pN\rbox[\pR]\mac{P}_{*,*} 
\quad \text{and} \quad
\mac{Q}_{*,*}\rbox[\pR]\mac{P}_{*,*}\lrarrow \mac{Q}_{*,*}\rbox[\pR]\pM,
\]
respectively. 
\end{thm}

In a similar fashion, by setting
\begin{equation}\label{biext}
\begin{split}
\mac{D}^{s,\tindex} &= \pimac[-s-\tindex]F(h^{s})\\
\mac{E}^{s,\tindex} &= \pimac[-s-\tindex]F(g^{s}),
\end{split}
\end{equation}
we obtain a cohomologically graded exact couple and maps 
of exact couples from both 
(\ref{leftext}) and~(\ref{rightext}) to~(\ref{biext}).  The relation
between these three exact couples is described by the following
result, which implies the remaining claims made in
Theorem~\ref{siso}. 

\begin{thm}\label{thmbiext}
The spectral sequence derived from~(\ref{biext}) has its $E_{1}$
complex canonically isomorphic to the total complex of
$\rfunc{\pR}{\mac{O}_{*,*}}{\mac{I}_{*}^{*}}$.  Moreover, under the canonical
isomorphisms, the maps into this complex from the $E_{1}$ complexes of (\ref{leftext}) and
(\ref{rightext}) coincide with the augmentations
\[
\rfunc{\pR}{\mac{O}_{*,*}}\pM\lrarrow
\rfunc{\pR}{\mac{O}_{*,*}}{\mac{I}_{*}^{*}} \quad \text{and} \quad
\rfunc{\pR}\pL{\mac{I}_{*}^{*}}\lrarrow
\rfunc{\pR}{\mac{O}_{*,*}}{\mac{I}_{*}^{*}},
\]
respectively. 
\end{thm}

For the proofs of these results, recall that, in the bar 
construction model of the homotopy colimit, $T_{n}$ is the 
geometric realization of a simplicial 
object which in simplicial degree $d$ is the wedge product, indexed on 
the set of $d$ composable
arrows
\[
(s_{0},t_{0})\from \cdots \from (s_{d},t_{d})
\]
in $\dcat{n}$, of the summands $\bar N_{s_{d}}\sma_{R}\bar M_{t_{d}}$.  The face
maps are given by dropping arrows on either end, composing arrows in
the middle, and by 
the action of $\dcat{n}$ on $\bar N_{s_{d}}\sma_{R}\bar M_{t_{d}}$.
The degeneracy maps simply insert identity maps. 
Likewise, the homotopy limit $U^{n}$ is the geometric realization (or
Tot) of a cosimplicial object which in cosimplicial degree $d$ is the
product, indexed on the set of $d$ composable arrows in $\dcat{n}$, of the factors 
$F_{R}(\bar L_{s_{d}},\bar M^{t_{d}})$.  Here, the face maps are given by
dropping or composing arrows and by the contravariant
action of $\dcat{n}$ on $F_{R}(\bar L_{s_{d}},\bar M^{t_{d}})$.  
The degeneracy maps again simply insert identity maps.

\begin{proof}[Proof of Theorems~\ref{thmbitor} and~\ref{thmbiext}]
We give the details of the proof of Theorem~\ref{thmbitor}.  The proof 
of Theorem~\ref{thmbiext} is quite similar.
The canonical maps from the homotopy colimits to the point-set colimits
induce a map $T_{0}\to\bar M_{0}\sma_{R}\bar N_{0}$
and, for $s> 0$, maps
\[
T_{s}/T_{s-1} \lrarrow (\bar N_{0})\sma_{R}(\bar M_{s}/\bar M_{s-1})
\vee \cdots \vee
(\bar N_{s}/\bar N_{s-1})\sma_{R}(\bar M_{0}).   
\]
These maps induce maps 
\begin{equation}\label{doublemap}
C(g_{s})
\lrarrow
\bigvee_{n+m=s} \dersma^{R}({C(\cofN{n})},{C(\cofM{m})})
\iso
\bigvee_{n+m=s} \dersma^{R}({\Sigma^{n}Q_{n}},{\Sigma^{m}P_{m}}).
\end{equation}
This gives a canonical map from $E^{1}$ to $\mac{Q}_{*,*}\rbox[\pR]
\mac{P}_{*,*}$. It is easy to see that this is a map of chain
complexes and that the maps 
\begin{gather*}
\pimac  C(g_{s})\lrarrow 
\pimac  C(\id_{N}\sma \cofM{s})\iso
\pN\rbox[\pR]\mac{P}_{s,*}, \\
\pimac  C(g_{s})\lrarrow
\pimac  C(\cofN{s}\sma \id_{M})\iso
\mac{Q}_{s,*}\rbox[\pR]\pM
\end{gather*}
factor through the augmentations
\[
\mac{Q}_{*,*}\rbox[\pR]\mac{P}_{*,*}\lrarrow \pN\rbox[\pR]\mac{P}_{*,*}
\quad \text{and} \quad 
\mac{Q}_{*,*}\rbox[\pR]\mac{P}_{*,*}\lrarrow \mac{Q}_{*,*}\rbox[\pR]\pM
\]
as required.

It remains to see that (\ref{doublemap}) induces an isomorphism on
homotopy groups.  Let $\dcat{n,m}$ denote the full subcategory of
$\dcat{n}$ whose objects are the pairs $(s,t)$ with $s\leq m$.  Also, let
$T_{n,m}=\Hocolim_{\dcat{n,m}}\bar N_{s}\sma_{R}\bar M_{t}$.  Standard
homotopy theory arguments show that the canonical map from the homotopy
colimit to the colimit induces weak equivalence of the cofiber of
$T_{s-1,0}\to T_{s,0}$ with $\bar N_{0}\sma \bar M_{s}$ and the
cofiber of the map $T_{s,m-1}\to T_{s,m}$ with $(\bar N_{m}/\bar
N_{m-1})\sma_{R}\bar M_{s-m}$.  A filtration argument now finishes
the proof.
\end{proof}


We close this section with an explanation of the 
Yoneda pairing of Hyper-Ext spectral sequences.  
Assume that $K$ is yet another left $R$-module, and let $\pK=\pimac K$.  
Also, let
$\mac{E}_{r}^{s,\tindex}(L,K)$ denote the spectral sequence for 
$\TExt^{*}_{R}(L,K)$ derived from~(\ref{leftext}), 
$\mac{E}_{r}^{s,\tindex}(K,M)$ denote the spectral sequence for 
$\TExt^{*}_{R}(K,M)$ derived
from~(\ref{rightext}), and $\mac{E}_{r}^{s,\tindex}(L,M)$ denote 
the spectral sequence for $\TExt^{*}_{R}(L,K)$ derived
from~(\ref{biext}). A map of complexes
\begin{equation}\label{erpair}
\bigoplus_{\ell+m=s}
\mac{E}_{r}^{m,\tindex}(K,M)\boxprod 
\mac{E}_{r}^{\ell,\tindex}(L,K)
\lrarrow\mac{E}_{r}^{s,\tindex}(L,M)
\end{equation}
induces in the usual way a map 
\begin{eqnarray*}
\lefteqn{\bigoplus_{\ell+m=s}
\mac{E}_{r+1}^{m,\tindex}(K,M)\boxprod 
\mac{E}_{r+1}^{\ell,\tindex}(L,K)}\hspace{6em}
\\
&=&
\bigoplus_{\ell+m=s}
H^{m}(\mac{E}_{r}^{*,\tindex}(K,M))\boxprod 
H^{\ell}(\mac{E}_{r}^{*,\tindex}(L,K))
\\
&\lrarrow&
H^{s}(\mac{E}_{r}^{*,\tindex}(K,M)\boxprod
\mac{E}_{r}^{*,\tindex}(L,K)) \\
&\lrarrow& H^{s}(\mac{E}_{r}^{*,\tindex}(L,M))
\\
&=& 
\mac{E}_{r+1}^{s,\tindex}(L,M).  
\end{eqnarray*}
of graded Mackey functors for $E_{r+1}$.  A pairing of spectral sequences is 
map of complexes~(\ref{erpair}) for
$E_{1}$ such that the induced map on $E_{2}$ is 
a map of complexes, and, inductively, the induced map on
each $E_{r+1}$ is a map of complexes.

\begin{thm}\label{ssypair}
The composition pairing
\[
\rfunc{\pR}{\pK}{\mac{I}_{*}^{m}}\boxprod 
\rfunc{\pR}{\mac{O}_{\ell,*}}{\pK}
\lrarrow \rfunc{\pR}{\mac{O}_{\ell,*}}{\mac{I}_{*}^{m}}
\]
induces a pairing of the Hyper-Ext spectral sequences.  On
$E_{2}$, this pairing agrees with the Yoneda pairing.  Also, this pairing agrees on 
$E_{\infty}$ with the associated graded of
the composition pairing for $\TExt_{R}^{*}$.
\end{thm}

\begin{proof}
The maps 
\begin{equation}\label{compmaps}
F_{R}(K,\bar M^{m})\sma F_{R}(\bar L_{\ell},K)\lrarrow
F(\bar L_{\ell},\bar M^{m}) \lrarrow U^{\ell+m}
\end{equation}
of $S$-modules, orthogonal spectra, or symmetric spectra induce maps 
\[
\derF_{R}({K},{\Omega^{m}I^{m}}) \sma_{S}
\derF_{R}({\Sigma^{\ell}O_{\ell}},{K}) \lrarrow F(g^{\ell+m})
\]
in the stable category.  On homotopy groups, these maps induce 
the pairing on $E_{1}$ described in the theorem, which on $E_{2}$ is
the Yoneda pairing, by definition.  The composition maps 
\eqref{compmaps} also induce maps 
\begin{gather*}
\derF_{R}({K},{F(\bar M^{m}\to \bar M^{m-r})}) \sma 
\derF_{R}({C(\bar L_{\ell-r}\to \bar L_{\ell})},{K})\hspace*{2cm}\\ 
\hspace*{5cm}\lrarrow F(U^{\ell+m}\to U^{\ell+m-r})
\end{gather*}
in the stable category relating the fiber and cofiber of the canonical
maps $\bar M^{m}\to \bar M^{m-r}$ and $\bar L_{\ell-r}\to \bar
L_{\ell}$,  
respectively, to the fiber of the canonical map $U^{\ell+m}\to U^{\ell+m-r}$. 
An inductive 
argument using these maps indicates that the pairing on $E^{r}$ preserves the
differential.  The statement about $E_{\infty}$ is clear.
\end{proof}

When $R$ is a commutative $S$-algebra, the pairing described 
in Theorem \ref{ssypair} descends to a pairing defined in terms 
of $\rbox[\pR]$ rather than $\boxprod$.  Moreover, the identification 
of the pairing on the $E_{\infty}$ level respects the $\pR$-module 
structures.  Theorem \ref{ssypair}, together with these observations, 
establishes the properties of the Yoneda pairing of spectral 
sequences noted in the introduction.

\appendix

\section{The smash product and the box product}\label{secpistar}

The purpose of this appendix is to prove Theorem~\ref{pisym}, which asserts that the
homotopy Mackey functor is lax symmetric monoidal.  Since the argument
in the case of \Z-graded Mackey functors is well-known, we concentrate
on the \ROG-graded case.   

As in Section~\ref{sectop}, choose a model $S^{\tindex}$ for the 
\tindex-sphere for each element \tindex of \ROG with the restrictions that
$S^{0}=\gsS$ and that $S^{\tindex}$ is the
smash product of $\gsS$ with the standard $n$-sphere space whenever \tindex 
is the trivial representation of dimension $n>0$.  For all
other \tindex, the object $S^{\tindex}$ may be chosen arbitrarily from
the appropriate homotopy class.  When $\tindex=\aindex +\bindex$, 
$S^{\aindex}\sma S^{\bindex}$ and $S^{\tindex}$ are isomorphic in the
equivariant stable category.  Thus, we can choose an isomorphism 
\[
f_{\aindex,\bindex}\colon S^{\aindex +\bindex}\lrarrow S^{\aindex}\sma S^{\bindex}.
\]
Such a choice gives a homomorphism
\begin{eqnarray*}
\shom{S^{\aindex}\sma X_{+}}{\gsM}\otimes
\shom{S^{\bindex}\sma Y_{+}}{\gsN}
&\lrarrow&
\shom{S^{\aindex}\sma X_{+}\sma S^{\bindex}\sma Y_{+}}{\gsM \sma \gsN}\\
&\lrarrow&
\shom{S^{\aindex+\bindex}\sma (X\times Y)_{+}}{\gsM \sma \gsN},
\end{eqnarray*}
which is natural in $\gsM$ and $\gsN$ in the equivariant stable category and
$X$ and $Y$ in the Burnside category.  The universal property of the
box product then gives a natural map
\[
\phi \colon \pimac\gsM \gbxprd \pimac\gsN\lrarrow \pimac(\gsM\sma \gsN).
\]
There is also a canonical natural map ${\iota}\colon {\mgMBr}\to{\pimac S}$ that
includes $\mMBr$ as $\pimac[0]S$.  This map does not require
any choices.  Theorem~\ref{pisym} is merely a less detailed 
version of the following result.  

\begin{thm}\label{apppisym}
The maps $f_{\aindex,\bindex}$ may be chosen so that $\phi$ together
with ${\iota}\colon {\mgMBr}\to{\pimac S}$ provide the functor $\pimac$ with a
lax symmetric monoidal structure.  
\end{thm}

In other words, the maps $f_{\aindex,\bindex}$ may be chosen so that the unit diagrams
\begin{equation}\label{piunit}
\begin{gathered}
\xymatrix{
\mgMBr \gbxprd \pimac \gsN\ar[d]_{\lambda}\ar[r]^{\id \gbxprd \iota}
&\pimac S\gbxprd \pimac \gsN\ar[d]^{\phi}
&\pimac \gsM \gbxprd \mgMBr \ar[d]_{\rho}\ar[r]^{\iota \gbxprd \id}
&\pimac\gsM \gbxprd \pimac S \ar[d]^{\phi}\\
\pimac\gsN&\pimac (S\sma \gsN)\ar[l]^{\pimac\lambda}
&\pimac\gsM&\pimac (\gsM \sma S)\ar[l]^{\pimac\rho},
}
\end{gathered}
\end{equation}
the associativity diagram
\begin{equation}\label{piassoc}
\begin{gathered}
\xymatrix@C=-5em{
&\pimac (\gsL \sma \gsM) \gbxprd \pimac\gsN
  \ar[rr]^{\phi}
&\hskip13em
&\pimac ((\gsL \sma \gsM) \sma \gsN)\ar[dr]^{\iso}\\
(\pimac \gsL \gbxprd \pimac\gsM)\gbxprd \pimac\gsN
  \ar[ur]^{\phi\gbxprd \id}\ar[dr]_{\iso}
&&&&\pimac (\gsL \sma (\gsM \sma \gsN)),\\
&\pimac \gsL \gbxprd(\pimac\gsM \gbxprd \pimac\gsN)
  \ar[rr]_{\id\gbxprd\phi}
&&\pimac \gsL \gbxprd \pimac(\gsM \sma \gsN) \ar[ur]_{\phi}
}
\end{gathered}
\end{equation}
and the symmetry diagram
\begin{equation}\label{picommut}
\begin{gathered}
\xymatrix{
\pimac \gsM \gbxprd \pimac \gsN\ar[d]_{\phi}\ar[r]^{\iso}
&\pimac \gsN \gbxprd \pimac \gsM\ar[d]^{\phi}\\
\pimac (\gsM \sma \gsN)\ar[r]_{\iso}
&\pimac (\gsN \sma \gsM)
}
\end{gathered}
\end{equation}
commute.  The redundant right unit diagram has been included 
because, with it, diagrams \eqref{piunit} and \eqref{piassoc} 
together imply that $\pimac$ is lax monoidal.  This is verified 
before the question of commutativity 
is addressed.  A direct argument in
terms of specifying how to choose the isomorphisms
$f_{\aindex,\bindex}$ appears to be possible, but would require 
checking a long list of complicated details.  Instead we take a more
abstract approach that reduces to checking the vanishing of a certain
characteristic cohomology class defined in \cite[\S7]{breen2gerbes}.  
The definition of that class in our specific case is reviewed below.

Since $S^{0}=S$ has been chosen as the unit for the smash product in
the equivariant stable category, the canonical choices for the maps
$f_{0,\tindex}$ and $f_{\tindex,0}$ are the inverses of the unit
isomorphisms $\lambda \colon S^0 \sma S^{\tindex} \to S^{\tindex}$ and
$\rho \colon S^{\tindex} \sma S^0\to S^{\tindex}$.  This is the unique 
choice making
the unit diagrams~(\ref{piunit}) commute.  This choice is assumed in  
the remainder of our argument.

It is easy to verify that the required associativity diagram~(\ref{piassoc})
commutes if and only if the diagram 
\begin{equation}\label{sphereassoc}
\begin{gathered}
\xymatrix@R=1.5em{ 
&{S^{\aindex + \bindex + \cindex} 
  \ar[dr]^{f_{\aindex + \bindex,\cindex}} 
  \ar[dl]_{f_{\aindex,\bindex + \cindex}} }\\
{S^{\aindex} \sma S^{\bindex + \cindex} 
  \ar[d]_{\id \sma f_{\bindex,\cindex}}} 
&&{S^{\aindex + \bindex} \sma S^{\cindex} 
  \ar[d]^{f_{\aindex,\bindex} \sma\, \id}}  \\
{S^{\aindex} \sma (S^{\bindex} \sma S^{\cindex})}\ar[rr]_{\iso}
&&{(S^{\aindex} \sma S^{\bindex}) \sma S^{\cindex}} }
\end{gathered}
\end{equation}
commutes for all $\aindex$,$\bindex$, and $\cindex$ in \ROG.
In general, the failure of this diagram to commute can be measured by
a unit in the Burnside ring.  To see this, recall that, in any symmetric 
monoidal additive category, the abelian group of
maps between any two objects is a module over the ring of
endomorphisms of the unit object (see for example the discussion of
``signs'' in Section~\ref{secgrdmackey}).  The abelian group 
\shom{S^{\aindex+\bindex+\cindex}}{S^{\aindex}\sma S^{\bindex}\sma
S^{\cindex}} is a one-dimensional free module over
$\shom{S}{S}=\mMBr(G/G)$ and is generated by any isomorphism.  The 
right-hand composite isomorphism in diagram~(\ref{sphereassoc}) 
is therefore the product of left-hand composite isomorphism and 
a well-defined unit $\cocycle{\aindex,\bindex,\cindex}$ in
$\mMBr(G/G)$.  Equivalently, $\cocycle{\aindex,\bindex,\cindex}$ can be 
defined as the unique unit in
$\mMBr(G/G)$ such that
\begin{equation}\label{assoccocycle}
\cocycle{\aindex,\bindex,\cindex}\cdot\id_{S^{\aindex +\bindex +\cindex}}
= f^{-1}_{\aindex,\bindex + \cindex}
\circ (\id \sma f_{\bindex,\cindex})^{-1} 
\circ a
\circ (f_{\aindex,\bindex} \sma\, \id) 
\circ f_{\aindex + \bindex,\cindex}.
\end{equation}
Here, $a$ denotes the associativity isomorphism in the equivariant
stable category.
Clearly, diagram~(\ref{sphereassoc}) commutes if and only if
$\cocycle{\aindex,\bindex,\cindex}=1$. Thus, we have proven:  

\begin{prop}
The associativity diagram~(\ref{piassoc}) commutes for all $L,M,N$ if and
only if $\cocycle{\aindex,\bindex,\cindex}=1$ for all
$\aindex,\bindex,\cindex\in \ROG$. 
\end{prop}

The isomorphism $\cocycle{\aindex,\bindex,\cindex}\cdot\id_{S^{\aindex
+\bindex +\cindex}}$ has a straight-forward structural interpretation
in terms of the choice of models $S^{\tindex}$ and isomorphisms
$f_{\aindex,\bindex}$.  Consider the full subcategory $\INV$ of the
equivariant stable category consisting of objects isomorphic to
spheres $S^{\tindex}$ (for $\tindex \in \ROG$), and let $\SPH$ be the
full subcategory consisting of the chosen models $S^{\tindex}$.  The
inclusion of $\SPH$ in $\INV$ is obviously an equivalence of categories.  The
smash product on the equivariant stable category 
restricts to a symmetric monoidal product on $\INV$.  The 
isomorphisms $f_{\aindex,\bindex}$ can be used 
to construct an equivalent symmetric monoidal smash product on 
$\SPH$.  There, the smash product of $S^{\aindex}$ and
$S^{\bindex}$ is $S^{\aindex+\bindex}$.  For maps
$g\colon S^{\aindex}\to S^{\aindex'}$ and $h\colon S^{\bindex}\to
S^{\bindex'}$, the induced map $g\sma h\colon S^{\aindex +\bindex}\to
S^{\aindex'+\bindex'}$ in $\SPH$ is the composite
\[
S^{\aindex +\bindex}\lrarrow[f_{\aindex,\bindex}]
S^{\aindex}\sma S^{\bindex}\lrarrow[g\sma h]
S^{\aindex'}\sma S^{\bindex'}\lrarrow[f^{-1}_{\aindex',\bindex'}]
S^{\aindex' +\bindex'}.
\]
Because of our restriction that $f_{0,\tindex}$ and $f_{\tindex,0}$ are the
inverses of unit isomorphisms, the unit isomorphisms in $\SPH$ are the 
appropriate identity maps.  The associativity isomorphisms in $\SPH$ are precisely
the maps $\cocycle{\aindex,\bindex,\cindex}\cdot
\id_{S^{\aindex+\bindex+\cindex}}$.  Since $\mMBr(G/G)$ is
commutative and composition in $\SPH$ is bilinear over $\mMBr(G/G)$, the
pentagon and triangle laws \cite[VII\S1]{CatWrk} for the monoidal 
structure on $\SPH$ translate into the following assertion about the elements
$\cocycle{\aindex,\bindex,\cindex}$:  

\begin{prop}
The elements $\cocycle{\aindex,\bindex,\cindex}$ satisfy the
cocycle condition 
\begin{equation*}
\cocycle{\aindex,\bindex,\cindex+\dindex} \cdot
\cocycle{\aindex+\bindex,\cindex,\dindex} =
\cocycle{\bindex,\cindex,\dindex} \cdot
\cocycle{\aindex,\bindex+\cindex,\dindex} \cdot
\cocycle{\aindex,\bindex,\cindex},
\end{equation*}
in the group of units of the Burnside ring and are normalized in the
sense that $\cocycle{\aindex,\bindex,\cindex}=1$ when any of
$\aindex$, $\bindex$, or $\cindex$ are $0$.
\end{prop}

This may be rephrased as the assertion that the collection 
$\{\cocycle{\aindex,\bindex,\cindex}\mid \aindex,\bindex,\cindex 
\in \ROG\}$ specifies a normalized $3$-cocycle
for the group cohomology of $\ROG$ with coefficients in group 
$\Bunits$ of units of the Burnside ring.  

For the proof of Theorem \ref{apppisym}, we must consider how this cocycle 
changes when a different collection of isomorphisms $f_{\aindex,\bindex}$ 
is selected.  Let 
${f'_{\aindex,\bindex}}\colon {S^{\aindex+\bindex}}\to{S^{\aindex} \sma 
S^{\bindex}}$ be another such collection satisfying our restriction 
on the maps $f'_{0,\tindex}$ and $f'_{\tindex,0}$.  The composite 
$f'_{\aindex,\bindex}\circ f_{\aindex,\bindex}^{-1}\colon 
S^{\aindex +\bindex}\to S^{\aindex+\bindex}$
is $\dclass{\aindex,\bindex}\cdot\id_{S^{\aindex +\bindex}}$ for a
well defined unit $\dclass{\aindex,\bindex}$ in the Burnside ring.
It is easy to verify that the cocycle $(\cocycle{\aindex,\bindex,\cindex}')$
associated to the isomorphisms $f'_{\aindex,\bindex}$ is given by 
\[
\cocycle{\aindex,\bindex,\cindex}' =
\dclass{\aindex,\bindex}\cdot
\dclass{\aindex+\bindex,\cindex}\cdot
\dclass{\aindex,\bindex+\cindex}^{-1}\cdot
\dclass{\bindex,\cindex}^{-1}\cdot
\cocycle{\aindex,\bindex,\cindex}.
\]
When the collection $(\dclass{\aindex,\bindex})$ is regarded as a $2$-cochain for
group cohomology, its boundary satisfies 
\[
(d\dclass{})_{\aindex,\bindex,\cindex} = \dclass{\aindex,\bindex}
\dclass{\aindex,\bindex+\cindex}^{-1}\dclass{\aindex+\bindex,\cindex}
\dclass{\bindex,\cindex}^{-1}.  
\]
Changing the choice of the maps
$f_{\aindex,\bindex}$ therefore changes the class 
$(\cocycle{\aindex,\bindex,\cindex})$ by a coboundary, from which 
it follows that $\cocycle{}$ determines a well-defined element of
$H^{3}(\ROG;\Bunits)$.  

Conversely, given any normalized $2$-cochain
$(\dclass{\aindex,\bindex})$ for $\ROG$ with
coefficients in $\Bunits$, the rule
$f'_{\aindex,\bindex}=\dclass{\aindex,\bindex}\cdot
f_{\aindex,\bindex}$ gives a collection of isomorphisms 
${f'_{\aindex,\bindex}}\colon {S^{\aindex+\bindex}}\to{S^{\aindex} \sma
S^{\bindex}}$ with $f'_{0,\tindex}$ and $f'_{\tindex,0}$ the inverses
of the unit isomorphisms, and with associated cocycle
$\cocycle{}'=d\dclass{}\cdot\cocycle{}$.  An easy computation indicates
that  
changing the models $S^{\tindex}$ does not alter the cohomology 
class represented by the collection 
$(\cocycle{\aindex,\bindex,\cindex})$.  These observations are 
summarized in the following result:  

\begin{prop}\label{cohomology}
There exists a cohomology class $\coclass$ in $H^{3}(\ROG;\Bunits)$ such that, 
for any fixed choice of models $S^{\tindex}$, the
association~(\ref{assoccocycle}) defines a surjection from
the set of collections of isomorphisms
${f_{\aindex,\bindex}}\colon {S^{\aindex+\bindex}}\to{S^{\aindex} \sma
S^{\bindex}}$ (with $f_{0,\tindex}$ and $f_{\tindex,0}$ the inverses of
the unit isomorphisms) onto the set
of normalized $3$-cocycles in the cohomology class $\coclass$.
\end{prop}

This result indicates, if the cohomology class $\coclass$ vanishes, 
then any collection of isomorphisms 
$f'_{\aindex,\bindex}$ can be adjusted via some normalized 
$2$-cochain into a collection $f_{\aindex,\bindex}$ making 
diagram \eqref{sphereassoc} commute.  This reduction of the 
proof of Theorem \ref{apppisym} to a question in cohomology 
opens the way for us to restrict our attention to actual, 
rather than virtual, real representations.  Let \posROG be 
the commutative monoid of isomorphism
classes of actual real representations of $G$.  It follows from 
\cite[4.1]{segalgamma} that the inclusion of \posROG into \ROG
induces a homotopy equivalence from $B(\posROG)$ to $B(\ROG)$.  
The associated restriction map
\begin{equation*}
H^3(\ROG;\Bunits) \lrarrow H^3(\posROG;\Bunits)
\end{equation*}
is therefore an isomorphism.  Thus, it suffices to show that the
image $\poscoclass$ of the cohomology class $\coclass$ of
Proposition~\ref{cohomology} is trivial in $H^3(\posROG;\Bunits)$.

The process just described for associating a normalized $3$-cocycle
in the cohomology of $\ROG$ to any 
collection of maps $f_{\aindex,\bindex}$ indexed on 
$\aindex,\bindex \in \ROG$ works equally well to associate a 
normalized $3$-cocycle in the cohomology of $\posROG$ to any 
collection of maps $f_{\aindex,\bindex}$ indexed on 
$\aindex,\bindex \in \posROG$.  Moreover the resulting class in the 
cohomology of $\posROG$ is clearly the restriction of the 
class in the cohomology of $\ROG$.  This allows us to turn the 
whole argument backwards:  we prove that the class $\poscoclass$ 
is trivial by showing
that we can choose model spheres $S^{\tindex}$ for $\tindex\in
\posROG$ and isomorphisms $f_{\aindex,\bindex}$ for $\aindex,\bindex
\in \posROG$ such that the diagrams~(\ref{sphereassoc}) commute.

Let $\rho_1$, $\rho_2$, \ldots, $\rho_r$ be an enumeration of the
irreducible real representations of $G$.  For each $\rho_i$, select an
associated sphere $S^{\rho_i}$.  Any nonzero $\tindex \in \posROG$ can be
written as a sum
\begin{equation*}
\tindex = n_{1}\rho_{1}+\cdots +n_{r}\rho_{r}
\end{equation*}
in which at least one of the $n_i$ is nonzero.  Let $S^\tindex$ be $S
\sma (S^{\rho_{1}})^{(n_{1})} \sma \cdots \sma
(S^{\rho_{r}})^{(n_{r})}$.  Here, $(S^{\rho_{i}})^{(n_{i})}$ denotes
the $n_{i}$-fold smash product of copies of $S^{\rho_{i}}$.  For each
pair $\aindex$, $\bindex$ in \posROG, select the map
${\bar{f}_{\aindex,\bindex}}\colon
{S^{\aindex+\bindex}}\to{S^{\aindex} \sma S^{\bindex}}$ to be that
which uses the associativity and commutativity isomorphisms of the
equivariant stable category to rearrange the spheres
$(S^{\rho_{i}})^{(n_{i})}$ appearing in $S^{\aindex+\bindex}$ into the
proper order for $S^{\aindex} \sma S^{\bindex}$ with as few
transpositions as possible.  Clearly this yields a collection of maps
which are appropriately unital and associative.  The associated
cocycle is then identically $1$, and the cohomology class
$\poscoclass$ is therefore trivial.  We have therefore proven the 
following result.

\begin{prop}
The cohomology class $\coclass$ of Proposition~\ref{cohomology} is
trivial.  Thus, we can choose the isomorphisms
$f_{\aindex,\bindex}$ for $\aindex,\bindex\in \ROG$ so that $\phi$
satisfies~(\ref{piunit}) and~(\ref{piassoc}).
\end{prop}

We still need to show that the natural transformation $\phi$ is
symmetric, that is, that it satisfies~(\ref{picommut}).  Let
$u(\aindex,\bindex)$ be the unique unit in the Burnside ring
such that the composite map
\begin{equation*}
S^{\aindex + \bindex} \lrarrow[f_{\aindex,\bindex}] S^{\aindex}\sma
S^{\bindex} \iso S^{\bindex}\sma S^{\aindex}
\lrarrow[f^{-1}_{\bindex,\aindex}] S^{\aindex + \bindex}
\end{equation*}
is $u(\aindex,\bindex)\cdot \id_{S^{\aindex +\bindex}}$.
By examining the definition of
$\phi$ and the symmetry isomorphism for the box product, it is easy to
see that $\phi$ is symmetric if and only if
$u(\aindex,\bindex)=\sigma(\aindex,\bindex)$.  Clearly $u$ is
antisymmetric; that is, $u(\aindex,\bindex) u(\bindex,\aindex)=1$.  
The commutativity of the diagram~(\ref{sphereassoc}) implies that 
$u$ is bilinear; that is, that $u(\aindex + \bindex, \cindex) =
u(\aindex,\cindex) u(\bindex,\cindex)$.  Thus, to see that our choice
of signs $\sigma$ in Section \ref{secgrdmackey} is consistent with our
choice of the maps $f_{\aindex,\bindex}$, it suffices to check
that $u(\rho,\rho)=\sigma(\rho,\rho)$ for each irreducible 
representation $\rho$ and that $u(\rho,\rho')=1$ whenever $\rho$ 
and $\rho'$ are distinct irreducible
representations.  This follows from the discussion
above and the observation that the twist map on $S^{\rho} \sma
S^{\rho}$ is homotopic to the map 
\begin{equation*}
\map{-1 \sma \id}{S^{\rho} \sma S^{\rho}}{S^{\rho} \sma S^{\rho}}.
\end{equation*}

\begin{rem}
We hinted in Section~\ref{secgrdmackey} that there is some flexibility
in the units $u(\aindex,\bindex)$ determined by the symmetry
isomorphism $S^{\aindex +\bindex}\to S^{\bindex +\aindex}$ and
therefore in the choice of $\sigma(\aindex,\bindex)$.  The analysis
leading to Proposition~\ref{cohomology} shows that, for a
fixed choice of model spheres $S^{\tindex}$, the set of collections of
isomorphisms $f_{\aindex,\bindex}$ (with $f_{0,\tindex}$ and
$f_{\tindex,0}$ the inverses of the unit isomorphisms) is a torsor for
the normalized $2$-cochains of $\ROG$ with coefficients in
$\Bunits$. Thus, the set of collections that
satisfy~(\ref{sphereassoc}) is a torsor for the normalized
$2$-cocycles.  If we adjust the maps $f_{\aindex,\bindex}$ 
by a $2$-cocycle $\dclass{}$, then the 
element $u(\aindex,\bindex)$ changes by a factor of
$\dclass{\bindex,\aindex}^{-1}\cdot \dclass{\aindex,\bindex}$.  It
follows that, for any collection of isomorphisms
satisfying~(\ref{sphereassoc}), $u(\rho,\rho)$ must be the unit
described above.  However, we can find a collection that takes on any
desired set of values for $u(\rho,\rho')$ for distinct irreducibles
$\rho,\rho'$, subject only to the restriction
$u(\rho,\rho')u(\rho',\rho)=1$.
\end{rem}


\end{document}